\documentclass[11pt,twoside]{article}
\usepackage{pb-diagram}
\usepackage{latexsym}
\input{amssym.def}
\input{amssym}
\newfont{\fra}{eufm10 scaled 1095}
\newfont{\Bb}{msbm10 scaled 1095}
\newfont{\Bbg}{msbm10 scaled 1280}
\newcommand\CC{{\mbox{\Bb C}}}
\newcommand\RR{{\mbox{\Bb R}}}

\newcommand\ZZ{{\mbox{\Bb Z}}}

\newcommand\fg{{\frak g}}
\newcommand\fh{{\frak h}}
\newcommand\fri{{\frak i}}
\newcommand\fj{{\frak j}}
\newcommand\fl{{\frak l}}
\newcommand\fm{{\frak m}}
\newcommand\fn{{\frak n}}
\newcommand\fa{{\frak a}}
\newcommand\fb{{\frak b}}
\newcommand\fd{{\frak d}}
\newcommand\fr{{\frak r}}
\newcommand\fs{{\frak s}}
\newcommand\fz{{\frak z}}

\newcommand\cH{{\cal H}}
\newcommand\cA{{\cal A}}

\newcommand{\fsl}{\mathop{{\frak s \frak l}}}
\newcommand\osc{\frak o \frak s \frak c}

\newcommand{\so}{\mathop{{\frak s \frak o}}}
\newcommand{\End}{\mathop{{\rm End}}}
\newcommand{\GL}{\mathop{{\it GL}}}
\newcommand{\Hom}{\mathop{{\rm Hom}}}
\newcommand{\Dera}{\mathop{{\rm Der_a}}}
\newcommand{\Id}{\mathop{{\rm Id}}}

\newcommand{\Ad}{\mathop{{\rm Ad}}}

\newcommand{\Ker}{\mathop{{\rm ker}}}
\newcommand{\im}{\mathop{{\rm im}}}
\newcommand{\Coker}{\mathop{{\rm coker}}}

\newcommand{\diag}{\mathop{{\rm diag}}}

\newcommand{\Span}{\mathop{{\rm span}}}
\newcommand{\mod}{\mathop{{\rm mod}}}
\newcommand\ip{\mbox{$\langle\cdot \,,\cdot \rangle$}}
\newcommand\ipa{{\langle\cdot \,,\cdot \rangle_\fa}}

\newcommand\proof{{\sl Proof. }}
\newcommand{\qed}{\hspace*{\fill}\hbox{$\Box$}\vspace{2ex}}
\newtheorem{theo}{Theorem}[section]
\newtheorem{pr}{Proposition}[section]
\newtheorem{de}{Definition}[section]
\newtheorem{ex}{Example}[section]
\newtheorem{re}{Remark}[section]
\newtheorem{co}{Corollary}[section]
\newtheorem{lm}{Lemma}[section]
\begin{document}
\title{Metric Lie Algebras with Maximal Isotropic Centre}
\author{Ines Kath and Martin Olbrich}
\maketitle
\begin{abstract}
\noindent
We investigate a certain class of solvable metric Lie algebras. For this
purpose a theory of twofold extensions associated
to an orthogonal representation of an abelian Lie algebra is developed. Among other things,
we obtain a classification scheme for indecomposable metric Lie algebras
with maximal isotropic centre and the classification of metric Lie algebras
of index $2$.

\end{abstract}
\tableofcontents
\section{Introduction}

In the present paper we study a certain class of metric Lie algebras in a systematic way. 
Here a metric Lie algebra is a finite-dimensional real Lie algebra equipped with an 
invariant non-degenerate symmetric bilinear form.
An isomorphism
of metric Lie algebras is by definition a Lie algebra isomorphism 
which is in addition an isometry with respect
to the given inner products. 

A metric Lie algebra $(\fg,\ip)$ is called decomposable if
it contains a proper ideal $\fri$ which is non-degenerate (i.e. $\ip_{|\fri\times\fri}$ is 
non-degenerate), and indecomposable otherwise. 
Then any metric Lie algebra
is the orthogonal direct sum of indecomposable metric Lie algebras.
Thus, in order to understand metric Lie algebras one has to understand
the indecomposable ones, and the natural (but certainly too ambitious)
task would be to obtain a classification of indecomposable metric Lie algebras
up to isomorphism.
   
Of course, in the case of definite inner products there is nothing
left to do. Any indecomposable metric Lie algebra with a definite inner
product is either compact simple or one-dimensional. However, the situation
changes dramatically if we consider indefinite inner products.
Here the simple and one-dimensional Lie algebras constitute only a small part of the set of 
all indecomposable metric Lie algebras. 

In \cite{MR85} A.~Medina and Ph.~Revoy proved a first structure result. They considered the 
following construction, which starts with a metric Lie algebra and produces metric Lie 
algebras of higher dimension. 

Let $(\fb,\ip_\fb)$  be a metric Lie algebra, $\fl$ a Lie algebra with invariant
(possibly degenerate) symmetric bilinear form $\ip_\fl$ and $\rho: \fl\rightarrow 
\Dera(\fb,\ip_\fb)$ a representation of $\fl$ by antisymmetric derivations on 
$(\fb,\ip_\fb)$. Let $\fd_\rho(\fb,\fl)$ be the extension
\begin{eqnarray}\label{schnupfen}
0\longrightarrow \fl^*\longrightarrow\fd_\rho(\fb,\fl)  \longrightarrow\fb\rtimes_\rho\fl
\longrightarrow 0
\end{eqnarray}

of the semi-direct sum of $\fb$ and $\fl$ defined by 
the coadjoint representation of $\fl^*$ on the second factor of $\fb\rtimes _\rho\fl$  and the 
cocycle $\beta\in Z^2(\fb\rtimes _\rho\fl,\fl^*)$ which is given by 
$\beta(X_1+L_1,X_2+L_2)=\langle X_1,\rho(\cdot)X_2\rangle_\fb$ for $X_i\in\fb$ and 
$L_i\in\fl$, $i=1,2$. Define an inner product $\ip$ on $\fd_\rho(\fb,\fl)$ by 
$$\langle \alpha_1+X_1+L_1,\alpha_2+X_2+L_2\rangle= \langle X_1, X_2\rangle_\fb + \langle 
L_1, L_2\rangle_\fl +\alpha_1(L_2)+\alpha_2(L_1)$$ for $\alpha_i\in\fl^*$, $X_i\in\fb$ and 
$L_i\in\fl$, $i=1,2$.
Then $(\fd,\ip)$ is a metric Lie algebra and it is called a double extension of $\fb$ by 
$\fl$.

What Medina and Revoy found in \cite{MR85} is that any non-simple non-abelian indecomposable 
metric Lie algebra $(\fd,\ip)$ is isomorphic to a double extension $\fd_\rho(\fb,\fl)$,  where 
$\fl$ is one-dimensional or simple. The idea of the proof is to choose a minimal abelian 
ideal $\fri\not=0$ in $\fd$ and to take $\fl=\fd/\fri^\perp$, $\fb=\fri^\perp/\fri$, 
$\fl^*=\fri$. See also \cite{OF} for a nice presentation of the result and its applications to 
conformal field theory.

This structure result is very suitable for the classification of non-simple indecomposable 
Lorentzian Lie algebras (i.e. metric Lie algebras of signature $(1,q)$). The centre $\fz$ of 
such a Lie algebra is non-trivial. Furthermore, it is isotropic because of indecomposability. 
Hence, in this case $\fri=\fz$ is a canonical choice for a (one-dimensional) minimal abelian 
ideal. In particular, the Lie algebra is canonically isomorphic to a double extension. This 
makes it easy to study decomposability and to solve the isomorphism problem. The 
classification result is due to A.~Medina (\cite{M85}, see also Section \ref{S5} for an exact 
formulation of the result). Originally this result was not a consequence of the structure 
theorem, but the structure theorem arose as a generalization of its proof.  

In \cite{MR85} Medina and Revoy applied the structure theorem also for non-simple 
indecomposable metric Lie algebras of signature $(2,q)$. This yields a description of such 
Lie algebras as double extensions of (the already classified) Lorentzian Lie algebras by 
one-dimensional ones. Furthermore, Medina and Revoy determined the (outer) antisymmetric 
derivations of a Lorentzian Lie algebra.  So we know in principle how a non-simple 
indecomposable metric Lie algebras of signature $(2,q)$ looks like. However, this is not a 
classification. On one hand a metric Lie algebra of signature $(2,q)$ which is the double 
extension of a Lorentzian Lie algebra by a one-dimensional one can be decomposable and Medina 
and Revoy did not study which of the obtained Lie algebras are really indecomposable. On the 
other hand it is not clear which of the obtained Lie algebras are isomorphic. Contrary to the 
Lorentzian case the description of a metric Lie algebra of signature $(2,q)$ as a double 
extension is not canonical. So one can obtain isomorphic Lie algebras of signature $(2,q)$ 
even if one starts with non-isomorphic Lorentzian Lie algebras.
Proposition \ref{PA} of the present paper is more adapted to this situation. It 
contains only canonical choices and allows a better study of indecomposability and isomorphy. 
In particular, it yields (together with Proposition \ref{rhino} and the results of Section 
\ref{S5}) a complete classification of indecomposable metric Lie algebras of signature 
$(2,q)$.

The classification result for signature $(2,q)$ was already announced in \cite{BK02} and it 
was one of our motivations for writing this paper to present a proof of this result in a 
quite systematic
framework.  Besides the classification for signature $(2,q)$ one 
can also find a 
classification of indecomposable metric Lie algebras of dimension $n\le6$ in \cite{BK02} 
although the paper does not concentrate on classification results. The emphasis is more on the 
construction of Lie 
groups with pseudo-Riemannian bi-invariant metric which admits parallel spinor fields using 
double extensions. One series of examples (called multiple extensions) corresponds to metric 
Lie algebras having a maximal isotropic centre and was constructed inductively by double 
extensions by a one-dimensional Lie algebra starting with an abelian Euclidean Lie algebra. 
The treatment of these examples suggested that there should be a better description of metric 
Lie algebras with maximal isotropic centre than the description by double extensions. This 
was a further motivation for writing the present paper.

From a geometric point of view, one is interested in the underlying semi-Riemannian symmetric 
space defined by a Lie group with bi-invariant metric.
Of course, Lie groups are very special symmetric spaces, and one would
like to understand all of them. On the other hand,
any symmetric space is a homogeneous space of a distinguished Lie group,
its transvection group. Let us discuss the relation between metric Lie
algebras and semi-Riemannian symmetric spaces in little more detail.

Up to local isometry, a semi-Riemannian symmetric space is given by a triple
$(\fg,\ip,\sigma)$, where $(\fg,\ip)$ is a metric Lie algebra and $\sigma$
is an involutive automorphism of $(\fg,\ip)$ such that  
$\fg=[\fm,\fm]\oplus \fm$, where $\fm$ is the $-1$-eigenspace of $\sigma$. By abuse of 
notation let us call such
a triple simply a symmetric space. The signature of a symmetric space is
by definition the signature of the restriction of $\ip$ to $\fm$.

This description already supports the philosophy, which might be considered
as the main motivation for a geometer to deal with metric Lie algebras,
that any good method of investigation
of metric Lie algebras should have a refinement which brings the automorphism
$\sigma$ into play, and then yields corresponding results for symmetric
spaces. Thus, the first step in the study of symmetric spaces should
be the study of metric Lie algebras.

Though symmetric spaces can be considered as more complicated objects
than metric Lie algebras, in many respects one knows more about symmetric spaces
than about metric Lie algebras.  For instance, there is no thorough treatment of metric Lie 
algebras
comparable to \cite{CP80}. Moreover, the classification
of indecomposable symmetric spaces of index at most $2$ is known for more
then 30 years \cite{CW70}, \cite{CP70} (see \cite{N02} for some
necessary corrections). While the classification result for Lorentzian Lie algebras which was 
discussed above is already considerably
younger, the corresponding result for metric Lie algebras of index $2$ has
not appeared in the literature before the announcement given in \cite{BK02}.
  
These circumstances are not only an accident. The study of metric
Lie algebras cannot be reduced to that of symmetric spaces.
The
underlying symmetric space of (a Lie group associated to) a metric
Lie algebra $(\fg,\ip)$ is the triple $(\tilde\fg,\ip_{\tilde\fg},\sigma)$,
where
$$ \tilde \fg=\{[(X,Y)]\in (\fg\oplus\fg)/\diag (\fz\cap \fg^\prime)\:|\: 
X+Y \in \fg^\prime\}\ ,$$
$\sigma([(X,Y)])=[(Y,X)]$, and $\ip_{\tilde\fg}$ is induced by the 
form $\ip\oplus\ip$    
on $\fg\oplus\fg$. Here $\fz$ denotes the centre of $\fg$. We see that it is impossible to recover 
$\fg$
from its underlying symmetric space in general. What can be recovered is the scalar
product $\ip$ and the Lie algebra structure on $\fg/(\fz\cap\fg^\prime)$, only.
Thus non-isomorphic metric Lie algebras may correspond to isomorphic
symmetric spaces. Moreover, an indecomposable metric Lie algebra may give
rise to a decomposable symmetric space.  

Quite recently we learned that Berard Bergery and his collaborators
in Nancy study symmetric
spaces (in up to now unpublished work) based on a construction which we call twofold extension. 
This 
construction can be applied to metric Lie algebras as well and then becomes a close relative of the 
double extension described
above. Roughly speaking, one simply replaces the semi-direct sum
$\fd_\rho(\fb,\fl)$ in (\ref{schnupfen}) by a quite arbitrary extension of $\fl$ by $\fb$. In 
addition,
we require $\fb$ to be abelian (therefore denoted by $\fa$ in the following). For this one has to 
pay the price that the
Lie algebra $\fl$ can be more general than one-dimensional or simple. 
We feel that this method has
the potential to give much more insight into the structure of general
metric Lie algebras than 
the double extension method. In particular, it can be shown that 
any indecomposable non-simple metric Lie algebra can be written as such a twofold extension
in a canonical way.  

In the present paper we test this method 
in the simplest case, namely for abelian $\fl$. The precise construction
of twofold extensions in this case is given in Proposition \ref{Pd}. It
turns out that a relatively broad class of metric Lie algebras is 
already covered by this simple construction, namely those metric Lie
algebras satisfying $\fg''\subset\fz$ (Proposition \ref{PA}). This class
contains indecomposable metric Lie algebras with non-trivial maximal 
isotropic centre (Theorem \ref{T}) 
and indecomposable non-simple metric Lie algebras of index at most $2$ (Proposition \ref{rhino}).

In Sections \ref{gris} and \ref{S4} we solve the crucial isomorphism
and indecomposability problems for our twofold extensions. The key
ingredient 
is the description of equivalence classes of twofold extensions associated
to a fixed orthogonal representation $\rho:\fl\rightarrow \so(\fa,\ip_\fa)$
by a
certain cohomomology set $\cH^2_C(\fl,\fa)$ (Corollary \ref{triangle}). 
This result is
completely analogous to the classical identification of equivalence classes
of extensions of a Lie algebra $\fl$ by an $\fl$-module $\fa$ with the
Lie algebra cohomology group $H^2(\fl,\fa)$.
Eventually, we can identify the set of isomorphism classes of indecomposable metric Lie algebras in 
our class with the union of orbit spaces of the
actions of certain linear groups on sets made of the ``indecomposable points'' of the cohomology 
sets 
$\cH^2_C(\fl,\fa)$ (Proposition \ref{Pwieder}).

How the classification scheme provided by Proposition \ref{Pwieder}
can be made much more explicit in the case of maximal isotropic centre
is discussed in Section \ref{S5}. In particular, we give a complete classification
of indecomposable metric Lie algebras having an at most $3$-dimensional maximal
isotropic centre in Theorem \ref{iso23}. Eventually, the classification of indecomposable metric 
Lie algebras of index $2$
is presented in Theorem \ref{ii}.

After we had obtained the results of Section \ref{gris} we realized
that a nonlinear cohomology theory for orthogonal representations
of general Lie algebras has been already developed by Grishkov in 
\cite{grishkov}. His results suggest that at least the results of Section
\ref{gris} should have a direct analogue for non-abelian $\fl$.
This and related topics might be the content of a forthcoming paper.      


\section{A structure theorem}
\label{S2}
We start this section by giving a construction method of certain metric Lie algebras.  
\begin{pr}\label{Pd}
    Let $(\rho, \fa)$ be an orthogonal representation of an abelian Lie 
algebra $\fl$ on
    the semi-Euclidean vector space $(\fa, \ip_{\fa})$. Furthermore, choose 
a 3-form
    $\gamma \in \bigwedge^{3}\fl^*$ and a cocycle $\alpha \in
    Z^{2}(\fl,\fa)$ satisfying
    \begin{equation}\label{EK}
       \langle \alpha (L_{1},L_{2}),\alpha (L_{3},L_{4})\rangle +
        \langle \alpha (L_{2},L_{3}),\alpha (L_{1},L_{4})\rangle +
	 \langle \alpha (L_{3},L_{1}),\alpha (L_{2},L_{4})\rangle
	  =0
    \end{equation}
    for all $L_{1},L_{2},L_{3},L_{4}\in \fl$.
    Then the bilinear map
    $$[\,\cdot \,,\cdot \,]:\ (\fl^{*}\oplus\fa\oplus\fl)^{2}\longrightarrow
    \fl^{*}\oplus\fa\oplus\fl$$
    defined by
    \begin{eqnarray*}
       \ \fl^{*}&\subset&\fz(\fl^{*}\oplus\fa\oplus\fl)\\
       \ [A_{1},A_{2}]&=&\langle \rho (\cdot)A_{1},A_{2}\rangle\in\fl^{*}\\
       \ [A,L]&=&\langle A,\alpha (L,\cdot\,)\rangle -L(A)\in\fl^{*}+\fa\\
       \ [L_{1},L_{2}]&=&\gamma(L_{1},L_{2},\cdot\,)+\alpha
        (L_{1},L_{2})\in\fl^{*}+\fa
    \end{eqnarray*}
    for all $L,L_{1},L_{2}\in\fl$, $A,A_{1},A_{2}\in\fa$ is a Lie bracket
    on $\fl^{*}\oplus\fa\oplus\fl$ and the bilinear form $\ip$ on
    $\fl^{*}\oplus\fa\oplus\fl$ defined by
    $$\langle Z_{1}+A_{1}+L_{1}, Z_{2}+A_{2}+L_{2}\rangle = \langle
    A_{1},A_{2}\rangle_{\fa}  +Z_{1}(L_{2}) + Z_{2}(L_{1})$$
    for all $Z_{1},Z_{2}\in\fl^{*}$, $A_{1},A_{2}\in\fa$ and
    $L_{1},L_{2}\in\fl$ is an ad-invariant (non-degenerate) inner
    product on $\fl^{*}\oplus\fa\oplus\fl$.
\end{pr}
\begin{de}  We denote the metric Lie algebra $(\,\fl^{*}\oplus\fa\oplus\fl,\
   \ip\,)$ constructed above by $\fd_{\alpha,\gamma}(\fa,\fl,\rho)$ and call it a twofold   
   extension.
\end{de}
Already in the case $\fa=0$ this construction leads to interesting metric
Lie algebras. For formal reasons we also include the trivial case $\fl=0$
into our considerations.\\

{\sl Proof of Proposition \ref{Pd}.} Consider the extension
$$0\longrightarrow \fa\longrightarrow \fd_{0}\longrightarrow 
\fl\longrightarrow
0$$
defined by $\rho$ and $\alpha \in Z^{2}(\fl,\fa)$.
Now define a 2-form
$\beta$ on $\fd_{0}=\fa\oplus\fl$ with values in $\fl^{*}$ by
    \begin{eqnarray*}
	\beta(A,L_{1})(L_{2})&=&\langle A,\alpha(L_{1},L_{2})\rangle_{\fa}\\
	\beta(A_{1},A_{2})(L)&=&\langle L(A_{1}),A_{2}\rangle_{\fa}\\
	\beta(L_{1},L_{2})(L_{3})&=&\gamma(L_{1},L_{2},L_{3})
    \end{eqnarray*}
for all $L,L_{1},L_{2},L_{3}\in\fl$ and $A,A_{1},A_{2}\in\fa$.
Using that $\alpha$ is a cocycle satisfying (\ref{EK}) it is not hard
to prove that $\beta$ is also a cocycle, i.e. $\beta\in
Z^{2}(\fd_{0},\fl^{*})$. Then the bilinear map $[\,\cdot \,,\cdot \,]$ is
exactly the Lie bracket on $\fl^{*}\oplus\fa\oplus\fl$ defined by the
central extension
$$0\longrightarrow \fl^{*}\longrightarrow \fl^{*}\oplus\fa\oplus\fl
\longrightarrow \fd_{0}\longrightarrow
0$$
associated with $\beta\in Z^{2}(\fd_{0},\fl^{*})$.
The invariance of $\ip$ follows from the orthogonality of $\rho$.
\qed

By construction, $\fl^{*}$ is an isotropic subspace of the
centre $\fz(\fd)$  of $\fd:=\fd_{\alpha,\gamma}(\fa,\fl,\rho)$.

\begin{de}\label{reg}
We call $\fd=\fd_{\alpha,\gamma}(\fa,\fl,\rho)$ regular if and only if $\fz(\fd)=
\fl^*$.
\end{de}

\begin{lm}\label{gaga}
$\fd_{\alpha,\gamma}(\fa,\fl,\rho)$ is not regular if and only if there
exist vectors $L_{0}\in\fl$ and $A_{0}\in \fa$, $L_0+A_0\ne 0$, such that
\begin{eqnarray}
    \rho(L_{0})&=&0\nonumber\\
    \rho(\,\cdot\,)(A_{0})&=&\alpha(L_{0},\cdot\,)\label{sing}\\
    \gamma(L_{0},\cdot\,,\cdot\,)&=&\langle A_{0},\alpha(\,\cdot \,,\cdot 
\,)\rangle\,. \nonumber
\end{eqnarray}
In other words, $\fd_{\alpha,\gamma}(\fa,\fl,\rho)$ is not regular if and only if Equations {\rm 
(\ref{sing})} are satisfied for some $L_0$, $A_0$ 
with $L_0\ne 0$ or the subspace $\fa^\fl\cap\alpha(\fl,\fl)^\perp\subset\fa$ is 
non-zero. Here $\fa^\fl$
denotes the space of invariants of $\rho$. 
\end{lm}
\proof
A straightforward computation shows that an element $X=L_0+A_0+Z_0\in 
\fd_{\alpha,\gamma}(\fa,\fl,\rho)=\fl\oplus \fa\oplus \fl^*$ is central if and only 
if the Equations (\ref{sing}) are satisfied.
\qed  

For any metric Lie algebra $\fg$ we have that 
$\fz(\fg)^\perp=\fg^\prime$.
Thus $\fd=\fd_{\alpha,\gamma}(\fa,\fl,\rho)$ is regular
if and only if $\fd^\prime=\fl^*\oplus\fa$.
In particular, in the regular case $\fz(\fd)\subset \fd'$, and $\fd'/\fz(\fd)$ is 
abelian.
Now we will prove a converse statement. Note that the centre of a non-abelian 
indecomposable metric Lie algebra $\fg$
is isotropic, i.e. $\fz(\fg)\subset \fg'$.

\begin{pr}\label{PA}
    If $(\fg,\ip)$ is a non-abelian indecomposable metric Lie algebra such that
    \begin{equation}\label{EA}
	\fg'/\fz(\fg) \mbox{ is abelian,}
    \end{equation}
    then there exist an abelian Lie
    algebra $\fl$, a semi-Euclidean vector
    space $(\fa, \ip_{\fa})$, an orthogonal representation $\rho$ of
    $\fl$ on $\fa$, a  cocycle $\alpha \in Z^{2}(\fl,\fa)$ satisfying 
(\ref{EK}),
    and a 3-form $\gamma \in
    \bigwedge^{3}\fl^*$ such that  $(\fg,\ip)$ is isomorphic to
    $\fd_{\alpha,\gamma}(\fa,\fl,\rho)$.

    In particular, we can choose
    \begin{eqnarray*}
	\fa&=& \fg'/\fz(\fg)\\
	\fl&=& \fg/\fg'
    \end{eqnarray*}
    and
    $$\rho(L)(A)=[l,a] + \fz(\fg)$$
    for $l\in\fg,\ a\in\fg',\ L=l+\fg',\
    A=a+\fz(\fg)$.
With this choice, $\fd_{\alpha,\gamma}(\fa,\fl,\rho)$ is regular.
\end{pr}
\proof
We first give the idea of the proof.
Choose $\fa,\, \fl$ and $\rho$ as above. Note, that because of (\ref{EA})
$\rho$ is correctly defined. We abbreviate the notation by $\fz=\fz(\fg)$.
By assumption the Lie algebra $\fg$ is the result of the following
two extensions with abelian kernel:
\begin{eqnarray*}
        &0\longrightarrow \fa\longrightarrow \fg/\fg'	\longrightarrow 
\fl\longrightarrow 0&\\
	&0\longrightarrow \fz\longrightarrow \fg\longrightarrow
	\fg/\fz	\longrightarrow 0.&
\end{eqnarray*}
The first extension is given by the representation $\rho$ as defined
above and, using an appropriate split, by a cocycle
$\alpha\in Z^{2}(\fl,\fa)$. The cocycle defining the second extension
will determine $\gamma \in \bigwedge^{3}\fl^*$.

Since $\fg'=\fz^\perp$ we have $\fa=
\fz^\perp/\fz$. Therefore $\ip$ induces a
non-degenerate inner product $\ip_{\fa}$ on $\fa$. Because of the
ad-invariance of $\ip$ the representation $\rho$ is orthogonal with
respect to $\ip_{\fa}$.

Because of $\fg'=\fz^\perp$ we can choose an isotropic complementary
subspace of $\fg'$ in $\fg$ and identify it with $\fl$ (as a vector
space). This yields a section
$$ s:\fl\longrightarrow\fg.$$
Since $s(\fl)$ is an isotropic complement of
$\fg'=\fz^{\perp}$ in $\fg$ the inner product $\ip$ defines a
dual pairing between $s(\fl)$ and $\fz$.
Therefore $s(\fl)^{\perp}\cap \fg'$ is a (non-degenerate) complement of
$\fz$ in $\fg'$. Identifying $\fa$ with
$s(\fl)^{\perp}\cap\fg'$ we obtain a section
$$t:\fa\longrightarrow \fg\ .$$

We define a cocycle $\alpha\in Z^{2}(\fl,\fa)$  by
$$\alpha(L_1,L_{2})=[s(L_{1}),s(L_{2})] + \fz. $$
Then $\alpha$ satisfies
Equation (\ref{EK}). This follows from
\begin{eqnarray*}
    \langle\alpha(L_1,L_2),\alpha(L_3,L_4)\rangle_{\fa}&=&
    \langle[s(L_{1}),s(L_{2})]\,,\,
    [s(L_{3}),s(L_{4})]\rangle\\&=&\langle[[s(L_{1}),s(L_{2})],
    s(L_{3})]\,,\,s(L_{4})]\rangle
\end{eqnarray*}
and the Jacobi identity in $\fg$. Finally we define a 3-form
$\gamma\in \bigwedge^{3}\fl^*$ by
$$\gamma(L_{1},L_{2},L_{3})=\langle[s(L_{1}),s(L_2)],s(L_3)\rangle.$$

Obviously we have an isomorphism
$$i:\fl^{*}\longrightarrow \fz$$
by $\langle i(Z),s(L)\rangle=Z(L)$ for $Z\in\fl^{*}$ and
$L\in\fl$. Now we will prove that the vector space isomorphism
$$i+t+s: \fl^{*}\oplus\fa\oplus\fl\longrightarrow \fg$$
is an isomorphism between the metric Lie algebras
$\fd_{\alpha,\gamma}(\fa,\fl,\rho)$ and $(\fg,\ip)$. It is almost
obvious that the map $i+t+s$ is an isometry. Thus it remains to show that it
is an isomorphism of Lie algebras. This follows from the following
equations which use the invariance of $\ip$:
\begin{eqnarray*}
    \ [A_{1},A_{2}]&=&\langle\rho(\cdot) A_{1},A_{2}\rangle_{\fa}\\
    \ [t(A_{1}),t(A_{2})]&\in&\fz\quad \mbox{ by (\ref{EA})}\\
    \ \langle [t(A_{1}),t(A_{2})],s(L)\rangle&=&\langle
    [s(L),t(A_{1})],t(A_{2})\rangle\,=\,\langle\rho(L)(A_{1}),A_{2}
    \rangle_{\fa}\\[1ex]
    \ [L,A]&=&-\langle A,\alpha(L,\cdot\,)\rangle_{\fa}+\rho(L)(A)\\
    \ [s(L),t(A)]&\in&\fg'\\
    \ [s(L),t(A)]&\equiv&t(\rho(L)(A))\mbox{ mod }\fz\\
    \langle [s(L),t(A)],s(L_{1})\rangle_{\fa}&=&-\langle
    t(A),[s(L),s(L_{1})]\rangle\,=\,-\langle
    A,\alpha(L,L_1)\rangle_{\fa}\\[1ex]
    \ [L_{1},L_{2}]&=&\gamma(L_{1},L_{2},\cdot\,)+\alpha(L_{1},L_{2})\\
    \ [s(L_{1}),s(L_{2})]&\in&\fg'\\
    \ [s(L_{1}),s(L_{2})]&\equiv&t(\alpha(L_{1},L_{2}))\mbox{ mod
    }\fz\\
    \ \langle
    [s(L_{1}),s(L_{2})],s(L_{3})\rangle&=&\gamma(L_{1},L_{2},L_{3}).
\end{eqnarray*}
\qed
\begin{re}{\rm
 It is natural to ask to what extend the tupel $(\alpha,\gamma)$
is determined by $\fg$. The above proof gives that $(\alpha,\gamma)$ depends
on the choice of a section $s:\fl\rightarrow \fg$. A closer view 
shows that $(\alpha,\gamma)$ depends on the induced section
$\bar s:\fl\rightarrow \fg/\fz$, only. This corresponds
to the fact that the orbit of $(\alpha,\gamma)$ under a certain action
of the group of 1-cocycles $C^1(\fl,\fa)$ is canonically attached to $\fg$ (see Section
\ref{gris}). 
}  
\end{re}

\begin{de}
    Let $\fri$ be an isotropic ideal in $(\fg,\ip)$. We will say that
    \begin{enumerate}
	\item $\fri$ is maximal isotropic if $\dim \fri = \min(p,q)$, where
	$(p,q)$ is the signature of $\ip$,
	\item $\fri$ is isomaximal if $\fri$ is not properly contained
	in a further isotropic ideal.
    \end{enumerate}
\end{de}
Obviously, each maximal isotropic ideal is isomaximal.
\begin{pr}\label{rhi}
Let $\fg$ be an indecomposable metric Lie algebra with non-trivial
isomaximal centre. Then $\fg$ satisfies condition {\rm (\ref{EA})}.
\end{pr}
\proof We have to show that
$\fz(\fg)^{\perp}/\fz(\fg)$ is abelian. In Lemma \ref{Lm} we will
prove even more, namely that $\fri^{\perp}/\fri$ is abelian for each
isomaximal ideal $\fri$.
\qed
\begin{lm}\label{L1}
    Let $(\fg,\ip)$ be an indecomposable metric Lie algebra and $\fr$
    the radical of $\fg$. Assume that $\fr\not=0$. Then
    $\fr^{\perp}\subset \fr$ holds.
\end{lm}
\proof Let $\fg=\fs\ltimes\fr$ be a Levi decomposition of $\fg$ and
consider $\fs_{1}:=\fs\cap\fr^{\perp}$. Clearly, $\fs_{1}$ is an ideal
in $\fs$. Since $\fs$ is semi-simple, there is a complementary
ideal $\fs_{2}$ of $\fs_{1}$ in $\fs$, i.e.
$\fs=\fs_{1}\oplus\fs_{2}$ as a direct sum of ideals. It is easy to
see that $\fs_{2}\ltimes\fr$ is a non-degenerate ideal in $\fg$.
Since $\fg$ is indecomposable and $\fr\not=0$ we obtain
$\fs_{2}\ltimes\fr=\fg$. Hence, $\fs_{1}=0$ and, consequently,
$\fr^{\perp}\subset\fr$.
\qed

\begin{lm} \label{Lm}
    Let $(\fg,\ip)$ be a metric Lie algebra with non-trivial
    radical $\fr\not=0$ and let $\fri\subset\fg$ be an isomaximal ideal. Then we have
    \begin{enumerate}
        \item $\fr\cap\fr^{\perp}\subset\fri$ .
    \end{enumerate}
    If, in addition,  $(\fg,\ip)$ is indecomposable, then
    \begin{itemize}
        \item[2.] $\fri^{\perp}\subset \fr$,
	\item[3.] $\fri^{\perp}/\fri$ is abelian.
    \end{itemize}
\end{lm}
\proof
Since the ideal $\fri$ is isotropic it is abelian. Hence, it is
contained in the radical $\fr$. In particular, $\langle
\fri,\fr\cap\fr^{\perp}\rangle=0$. Therefore, $\fri +
\fr\cap\fr^{\perp}$ is an isotropic ideal, which must be equal to
$\fri$ since $\fri$ is isomaximal. This implies
$\fr\cap\fr^{\perp}\subset \fri$.
Using this we obtain
$$\fri^{\perp}\subset
(\fr\cap\fr^{\perp})^{\perp}=\fr^{\perp}+\fr\subset\fr,$$
where the last inclusion uses Lemma \ref{L1}.
Now we are going to prove the third assertion. Since $\fri^\perp\subset \fr$
the ideal $\fri^{\perp}$ is solvable. Therefore, also
$$\fa:=\fri^{\perp}/\fri$$
is solvable. We define an ideal $\tilde\fz$ of $\fri^{\perp}$ by
$$\fri\subset\tilde\fz:=\{X\in\fri^{\perp}\mid [X,Y]\in\fri\mbox{ for
all } Y\in\fri^\perp\}\subset\fri^{\perp}.$$
Obviously,
$\tilde\fz/\fri$ is the centre of $\fa$. Consider now
$$ \fj:=[\fri^{\perp},\fri^{\perp}]\cap\tilde\fz+\fri.$$
Since $\fri, [\fri^{\perp},\fri^{\perp}]$ and $\tilde\fz$ are ideals
$\fj$ is also an ideal. We claim that $\fj$ is isotropic. Indeed,
$[\fri^{\perp},\fri^{\perp}]\cap\tilde\fz$ is isotropic, because of
$[\fri^{\perp},\fri^{\perp}]\perp\tilde\fz$. On the other hand,
$\fri$ is isotropic and $[\fri^{\perp},\fri^{\perp}]\cap\tilde\fz$ is
contained in $\fri^{\perp}$. Thus, $\fj$ is an isotropic ideal
containing $\fri$. Therefore, $\fj=\fri$, which implies
$[\fri^{\perp},\fri^{\perp}]\cap\tilde\fz\subset\fri$. Factorizing
this by $\fri$ yields
\begin{equation}\label{Eaa}
    [\fa,\fa]\cap\fz(\fa)=0,
\end{equation}
where $\fz(\fa)$ denotes the centre of $\fa$. On the other hand we have
\begin{equation}\label{Eza}
    [\fa,\fa]^{\perp}=\fz(\fa),
\end{equation}
Since $\ip$ induces a non-degenerate ad-invariant inner product on
$\fa$. Combining Equations (\ref{Eaa}) and (\ref{Eza}) we obtain
$$\fa=\fa'\oplus\fz(\fa),$$
which implies
$$\fa'=[\fa',\fa']=\fa^{(2)}.$$
Now we use that $\fa$ is solvable and obtain $\fa'=0$, which proves
the assertion.
\qed

In particular, we can apply Propositions \ref{PA} and \ref{rhi} to
indecomposable metric Lie algebras with maximal isotropic centre $\fz$. In this case the induced 
inner product on
$\fz^{\perp}/\fz$ is definite.  We may assume, perhaps after a sign change of $\ip$, that it is 
Euclidean. We obtain the following structure theorem.

\begin{theo}\label{T}
    If $(\fg,\ip)$ is an indecomposable metric Lie algebra of signature $(p,q)$, $p\le 
    q$, 
    with
    non-trivial maximal iso\-tropic centre, then there exist an abelian Lie
    algebra $\fl$, a Euclidean vector
    space $(\fa, \ip_{\fa})$, an orthogonal representation $\rho$ of
    $\fl$ on $\fa$, a  cocycle $\alpha \in Z^{2}(\fl,\fa)$ satisfying 
    {\rm(\ref{EK})},
    and a 3-form $\gamma \in
    \bigwedge^{3}\fl^*$ such that $\fd_{\alpha,\gamma}(\fa,\fl,\rho)$
is regular and $(\fg,\ip)$ is isomorphic to
    $\fd_{\alpha,\gamma}(\fa,\fl,\rho)$.
\end{theo}
\begin{pr}\label{rhino}
Let $\fg$ be a non-simple indecomposable metric Lie algebra
of index at most $2$. Then $\fg$ satisfies condition {\rm (\ref{EA})}.
\end{pr}

\begin{lm}\label{zero}
Let $\fg$ be as in Proposition \ref{rhino}. Then $\fg$ is solvable.
\end{lm}

\proof
First we observe that $\fg$ cannot be semi-simple.
Thus $\fg$ has a Levi
decomposition
$\fg=\fs\ltimes\fr$ with $\fr\ne 0$.
By Lemma \ref{L1} we have $\fr^\perp\subset \fr$. Therefore the
scalar product on $\fg$ induces a non-degenerate pairing between
$\fs$ and $\fr^\perp$. It follows that the index of $\fg$ is greater
or equal to $\dim \fs=\dim \fr^\perp$ which is at least $3$ if not zero.
We conclude that $\fs=0$.
\qed

\begin{lm}\label{un}
A nilpotent metric Lie algebra of index at most $1$ is abelian.
\end{lm}
\proof
This is known according to \cite{M85}, see also \cite{BK02}. For the
convenience of the reader we give a selfcontained argument in the
language developed in the present paper.
Let $\fn$ be a nilpotent metric Lie algebra of index at most $1$.
Then $\fn$ is the direct sum of a Euclidean abelian Lie algebra $\fn_0$
and, possibly, an indecomposable nilpotent Lie algebra $\fn_1$ of index 1.
We have to show that $\fn_1$ is abelian. Let $\fz$ be the centre of
$\fn_1$. It is one-dimensional and maximal isotropic in $\fn_1$.
By Theorem \ref{T}
we have $\fn_1\cong \fd_{\alpha,\gamma}(\fa,\fl,\rho)$ with
$\fl$ one-dimensional, $\alpha=\gamma=0$, and $\fa \cong \fn_1^\prime/\fz$
Euclidean. Let $H$ be a generator of $\fl$. Since $\fn_1$ is nilpotent the
anti-symmetric
endomorphism $\rho(H)$ has to be nilpotent. It follows that $\rho(H)=0$.
Therefore $\fn_1\cong \fd_{\alpha,\gamma}(\fa,\fl,\rho)$ is abelian.
\qed

{\sl Proof of Proposition \ref{rhino}. }
According to Lemma \ref{zero} the Lie algebra $\fg$ is solvable. By Lie's
Theorem $\fg^\prime$ is nilpotent. Thus $\fg^\prime/\fz$ is
a nilpotent metric Lie algebra of index at most $1$. It is abelian by Lemma
\ref{un}.
\qed

Combining Proposition \ref{rhino} with Proposition \ref{PA} we obtain
\begin{co}\label{i2}
If $(\fg,\ip)$ is a non-simple indecomposable metric Lie algebra of index at
most $2$, then there exist an abelian Lie
    algebra $\fl$,  $\dim\fl=\dim\fz$, a semi-Euclidean vector
    space $(\fa, \ip_{\fa})$, an orthogonal representation $\rho$ of
    $\fl$ on $\fa$, and a cocycle $\alpha \in Z^{2}(\fl,\fa)$ satisfying {\rm 
    (\ref{EK})},
    such that $\fd_{\alpha,0}(\fa,\fl,\rho)$ is regular and $(\fg,\ip)$ is 
    isomorphic to
    $\fd_{\alpha,0}(\fa,\fl,\rho)$.
\end{co}


\section{Equivalence of twofold extensions and the isomorphism problem}
\label{gris}

In order to approach the question whether two metric Lie algebras
$\fd_{\alpha_i,\gamma_i}(\fa_i,\fl_i,\rho_i)$, $i=1,2$, are isomorphic we first 
introduce a stronger equivalence relation on 
these objects. 

We consider extensions of a Lie algebra $\fh$ by an abelian
Lie algebra $\fa$
$$ 0\rightarrow \fa\stackrel{i}{\longrightarrow}\fg
\stackrel{p}{\longrightarrow} \fh\rightarrow 0\ .$$
If $\fh$ and $\fa$ are understood we denote such an extension by
$(\fg,i,p)$. Two such extensions $(\fg_j,i_j,p_j)$, $j=1,2$, are called
to be equivalent if and only if there 
exists a Lie algebra isomorphism $\Phi: \fg_1\rightarrow \fg_2$ such that
$\Phi\circ i_1=i_2$ and $p_2\circ\Phi=p_1$.  
It is well-known that $(\fg_1,i_1,p_1)\cong (\fg_2,i_2,p_2)$ if and only
if the corresponding representations of $\fh$ on $\fa$ coincide and
the cocycles $\alpha_j\in Z^2(\fh,\fa)$ defined after a choice of sections of $p_j$ 
are cohomologous. In particular, equivalence classes of extensions
with a fixed representation of $\fh$ on $\fa$ are in bijective correspondence
with elements of $H^2(\fh,\fa)$.   

Recall from the proof of Proposition \ref{Pd} that
$\fd_{\alpha,\gamma}:=\fd_{\alpha,\gamma}(\fa,\fl,\rho)$ is the result of two 
subsequent
extensions of Lie algebras
\begin{eqnarray*}
(\fd_\alpha,i,p):&& 0\rightarrow \fa\stackrel{i}{\longrightarrow}\fd_\alpha
\stackrel{p}{\longrightarrow} \fl\rightarrow 0\ ,\\
(\fd_{\alpha,\gamma},j,q):&&
0\rightarrow \fl^*\stackrel{j}{\longrightarrow}\fd_{\alpha,\gamma}
\stackrel{q}{\longrightarrow} \fd_\alpha\rightarrow 0\ ,
\end{eqnarray*}
where $\fl$ is abelian, and the representations of $\fl$ on $\fa$ and of
$\fd_\alpha$ on $\fl^*$ are orthogonal and trivial, respectively.

\begin{de}\label{dreist}
Let $\alpha_i\in Z^2(\fl,\fa)$, $i=1,2$, be cocycles satisfying {\rm (\ref{EK})},
and let $\gamma_i\in \bigwedge^{3}\fl^*$.
We call the metric Lie algebras $\fd_{\alpha_i,\gamma_i}(\fa,\fl,\rho)$
extension equivalent if and only if there exist
an equivalence of extensions  
$$\Phi: (\fd_{\alpha_1},i_1,p_1)\rightarrow
(\fd_{\alpha_2},i_2,p_2)$$  
and an isomorphism of metric Lie algebras
$$\Psi: \fd_{\alpha_1,\gamma_1}(\fa,\fl,\rho)\rightarrow
\fd_{\alpha_2,\gamma_2}(\fa,\fl,\rho)$$
making the following diagram commutative
\begin{equation}\label{bubble}
\dgARROWLENGTH=2em
\begin{diagram}
 \node{0}\arrow[2]{e}
 \node[2]{\fl^*}\arrow[3]{e,t}{j_1}\arrow[2]{s,r}{\Id}
 \node[3]{\fd_{\alpha_1,\gamma_1}}\arrow[2]{e,t}{q_1}\arrow[2]{s,r}{\Psi}
 \node[3]{\fd_{\alpha_1}}\arrow[2]{e}\arrow[2]{s,r}{\Phi}
 \node[2]{0}   
 \\ \\ 
 \node{0}\arrow[2]{e}
 \node[2]{\fl^*}\arrow[3]{e,t}{j_2}
 \node[3]{\fd_{\alpha_2,\gamma_2}}\arrow[3]{e,t}{q_2}
 \node[3]{\fd_{\alpha_2}}\arrow[2]{e}
 \node[2]{0}
\end{diagram}.  
\end{equation}
\end{de}

We now consider the standard Lie algebra cochain complex of the abelian
Lie algebra $\fl$ with values in $\fa$
$$0\rightarrow C^0(\fl,\fa)\stackrel {d}{\longrightarrow}C^1(\fl,\fa)\stackrel 
{d}{\longrightarrow}C^2(\fl,\fa)\stackrel {d}{\longrightarrow}\dots \stackrel 
{d}{\longrightarrow}C^l(\fl,\fa)\rightarrow 0\ ,$$
where $C^p(\fl,\fa)=\Hom(\bigwedge^p\fl,\fa)$, $l=\dim \fl$, and 
for $\tau\in C^p(\fl,\fa)$
$$ d\tau(L_1,\dots,L_{p+1})=\sum_{i=1}^{p+1} (-1)^{i-1} \rho(L_i) 
\tau(L_1,\dots,\hat L_i,\dots,L_{p+1}) \ .$$
The composition of maps
$$ C^p(\fl,\fa)\times C^q(\fl,\fa)\stackrel{\wedge}{\longrightarrow}
C^{p+q}(\fl,\fa\otimes\fa)\stackrel{\ip_\fa}{\longrightarrow} 
{\textstyle\bigwedge^{p+q}}\fl^*$$
defines a bilinear multiplication
$$ C^p(\fl,\fa)\times C^q(\fl,\fa) \rightarrow {\textstyle\bigwedge^{p+q}}\fl^*$$
which we will denote by $(\alpha,\tau)\mapsto \langle \alpha\wedge\tau\rangle_\fa$.
We are particularly interested in the case $p=2$, $q=1$. Then
$$ \langle \alpha\wedge\tau\rangle_\fa (L_1,L_2,L_3)=
\langle \alpha(L_1,L_2),\tau(L_3)\rangle_{\fa}
+\langle \alpha(L_3,L_1),\tau(L_2)\rangle_{\fa}
+\langle \alpha(L_2,L_3),\tau(L_1)\rangle_{\fa}\ .$$

\begin{pr}\label{abruest}
The metric Lie algebras $\fd_{\alpha_i,\gamma_i}(\fa,\fl,\rho)$ are
extension equivalent if and only if there exists a cochain
$\tau\in C^1(\fl,\fa)$ such that
\begin{equation}\label{walt}
\alpha_2-\alpha_1= d\tau
\end{equation}
and
\begin{equation}
\gamma_2-\gamma_1=\langle 
(\alpha_1+\textstyle\frac{1}{2}d\tau)\wedge\tau\rangle_{\fa}
\ .\label{disney}
\end{equation}
\end{pr}

\proof Assume that $\fd_{\alpha_i,\gamma_i}(\fa,\fl,\rho)$ are
extension equivalent, i.e., there exist isomorphisms $\Phi$ and $\Psi$ as in 
Definition \ref{dreist}. We use the vector space decomposition
$\fd_{\alpha_1}=\fl\oplus\fa=\fd_{\alpha_2}$ in order to define $\tau\in 
C^1(\fl,\fa)$ by
$\tau(L):=L-\Phi(L)$, $L\in\fl$. Because of $p_2(L-\Phi(L))=L-p_1(L)=0$ we
have $\tau(L)\in\fa$. Moreover, (\ref{walt}) holds. Indeed,
$$ \Phi^{-1}(L)=L+\tau(L)$$
and
\begin{eqnarray*}
\alpha_2(L_1,L_2)=[L_1,L_2]_{\fd_{\alpha_2}}&=&
[\Phi^{-1}(L_1),\Phi^{-1}(L_2)]_{\fd_{\alpha_1}}\\
&=&[L_1+\tau(L_1),L_2+\tau(L_2)]_{\fd_{\alpha_1}}
=\alpha_1(L_1,L_2)+d\tau(L_1,L_2)\ .
\end{eqnarray*}
We denote the commutator in $\fd_{\alpha_i,\gamma_i}(\fa,\fl,\rho)$ by
$[\cdot\,,\cdot]_i$. The decomposition 
$\fd_{\alpha_1,\gamma_1}(\fa,\fl,\rho)=\fl\oplus\fa\oplus\fl^*=
\fd_{\alpha_2,\gamma_2}(\fa,\fl,\rho)$ will us allow to consider $\fl\oplus
\fa$ as a subspace of $\fd_{\alpha_1,\gamma_1}(\fa,\fl,\rho)$ and of 
$\fd_{\alpha_2,\gamma_2}(\fa,\fl,\rho)$, 
$\Phi^{-1}$ as a map between $\fl\oplus
\fa$ and $\fd_{\alpha_1,\gamma_1}(\fa,\fl,\rho)$, etc.
We compute 
\begin{eqnarray}
\gamma_2(L_1,L_2,L)&=&\langle [L_1,L_2]_2, L\rangle\nonumber\\
&=&\langle [\Psi^{-1}(L_1),\Psi^{-1}(L_2)]_1, \Psi^{-1}(L)\rangle\nonumber\\
&=&\langle [\Phi^{-1}(L_1),\Phi^{-1}(L_2)]_1, \Phi^{-1}(L)\rangle\label{plumbum}\\
&=& \langle [L_1+\tau(L_1),L_2+\tau(L_2)]_1, L+\tau(L)\rangle\nonumber\\
&=& \langle [L_1,L_2]_1, L\rangle\nonumber\\
&&+ \langle [L_1,L_2]_1, \tau(L)\rangle
+ \langle [L_1,\tau(L_2)]_1, L\rangle
+\langle [\tau(L_1),L_2]_1, L\rangle\nonumber\\
&& +\langle [\tau(L_1),L_2]_1, \tau(L)\rangle
+ \langle [L_1,\tau(L_2)]_1, \tau(L)\rangle
+ \langle [\tau(L_1),\tau(L_2)]_1, L\rangle\nonumber\\
&=& \gamma_1(L_1,L_2,L)\nonumber\\
&&+\langle \alpha_1(L_1,L_2),\tau(L)\rangle_{\fa}
+\langle \alpha_1(L,L_1),\tau(L_2)\rangle_{\fa}
+\langle \alpha_1(L_2,L),\tau(L_1)\rangle_{\fa}\nonumber\\
&&+\langle \tau(L_{1}),L_{2}\tau(L)\rangle_{\fa}
+\langle \tau(L),L_{1}\tau(L_{2})\rangle_{\fa}
+\langle \tau(L_{2}),L\tau(L_{1})\rangle_{\fa}\nonumber\\
&=& \gamma_1(L_1,L_2,L)+ \langle 
\alpha_1\wedge\tau\rangle_{\fa}(L_1,L_2,L)+\textstyle\frac{1}{2}
\langle d\tau\wedge\tau\rangle_{\fa}(L_1,L_2,L)\ .\nonumber 
\end{eqnarray}
This proves (\ref{disney}).

Vice versa, let us assume that $\tau\in C^1(\fl,\fa)$ satisfying (\ref{walt})
and (\ref{disney}) is given. Let $\tau^*: \fa\rightarrow \fl^*$ be its adjoint,
where we have identified $\fa$ with $\fa^*$ using
$\ip_\fa$. With respect to the decompositions
$\fd_{\alpha_i}=\fl\oplus\fa$, $\fd_{\alpha_i,\gamma_i}=\fl\oplus\fa\oplus\fl^*$
we define
\begin{eqnarray}
\Phi(L+A)&:=&L-\tau(L)+A\nonumber\\
\label{Psitau}
\Psi(L+A+Z)&:=&\Phi(L+A)-\frac{1}{2}\tau^*\tau(L)+\tau^*(A)+Z\ .
\end{eqnarray}
Then we have
\begin{equation}\label{tsch} 
      \Psi = \left(\begin{array}{ccc}
                       \Id&\tau^*&-\frac12\tau^*\tau\\
                       0&\Id&-\tau\\
                       0&0&\Id
                    \end{array}  
              \right)         
\end{equation}
with respect to the decomposition $\fl^*\oplus\fa\oplus\fl$.
Because of (\ref{walt}) the map $\Phi$ defines an equivalence between 
$(\fd_{\alpha_1},i_1,p_1)$ and $(\fd_{\alpha_2},i_2,p_2)$.
Obviously, $\Psi$ satisfies
the relations (\ref{bubble}). Moreover, it is easy to verify that $\Psi$ is an 
isometry. 

It remains to show that
$\Psi$ is in fact a Lie algebra homomorphism. Using that $\Psi$ is an isometry
this is equivalent to
\begin{equation}\label{gum}
\langle [\Psi^{-1}(X_1),\Psi^{-1}(X_2)]_1, \Psi^{-1}(X)\rangle = 
\langle [X_1,X_2]_2, X\rangle
\end{equation}
for all $X_i=L_i+A_i+Z_i$, $X=L+A+Z\in \fl\oplus\fa\oplus\fl^*$.
We can assume that $Z_i=0$ since $Z_i$ is central and fixed by $\Psi^{-1}$
and that $Z=0$ since $[\Psi^{-1}(X_1),\Psi^{-1}(X_2)]_1$ is orthogonal to $\fl^*$. 
Observe that $\Psi^{-1}-\Phi^{-1}$ maps $\fl\oplus\fa$ to $\fl^*$. Using again that 
commutators are orthogonal to $\fl^*$, we obtain  
$$ \langle [\Psi^{-1}(X_1),\Psi^{-1}(X_2)]_1, \Psi^{-1}(X)\rangle = 
\langle [\Phi^{-1}(X_1),\Phi^{-1}(X_2)]_1, \Phi^{-1}(X)\rangle\ . $$
If $X=A\in\fa$, then this is equal to
\begin{eqnarray*} 
\langle [\Phi^{-1}(X_1),\Phi^{-1}(X_2)]_1, A\rangle 
&=& \langle [\Phi^{-1}(X_1),\Phi^{-1}(X_2)]_{\fd_{\alpha_1}}, A\rangle \\
&=&\langle \Phi^{-1}([X_1,X_2]_{\fd_{\alpha_2}}), A\rangle = 
\langle [X_1,X_2]_2, A\rangle\ ,
\end{eqnarray*} 
since $\Phi^{-1}$ is a Lie algebra homomorphism which is equal to the identity when 
restricted to $\fa$. 
 
Now let $X=L\in\fl$.
Employing that $\Phi^{-1}(L)=L+\tau(L)$ we compute
$$\langle [\Phi^{-1}(A_1),\Phi^{-1}(A_2)]_1, \Phi^{-1}(L)\rangle =
\langle [A_1,A_2]_1, L+\tau(L)\rangle =\langle A_2, LA_1\rangle_{\fa} =\langle 
[A_1,A_2]_2, L\rangle$$
and 
\begin{eqnarray*} 
\langle [\Phi^{-1}(L_1),\Phi^{-1}(A_2)]_1, \Phi^{-1}(L)\rangle &=&
\langle [L_1+\tau(L_1),A_2]_1, L+\tau(L)\rangle\\
&=&\langle L_1A_2, \tau(L)\rangle_{\fa}
+\langle \alpha_1(L,L_1), A_2\rangle_{\fa} -\langle \tau(L_1), LA_2\rangle_{\fa}\\
&=& \langle (\alpha_1+d\tau)(L,L_1),A_2\rangle_{\fa} = 
\langle \alpha_2(L,L_1),A_2\rangle_\fa\\
&=&\langle [L_1,A_2]_2, L\rangle\ .
\end{eqnarray*}
Eventually we have (repeating the computation (\ref{plumbum}))
\begin{eqnarray*} 
\langle [\Phi^{-1}(L_1),\Phi^{-1}(L_2)]_1, \Phi^{-1}(L)\rangle
&=& \langle [L_1+\tau(L_1),L_2+\tau(L_2)]_1, L+\tau(L)\rangle \\
&=&\gamma_1(L_1,L_2,L)+\langle (\alpha_1+\textstyle\frac{1}{2}
 d\tau)\wedge\tau\rangle_{\fa}(L_1,L_2,L)\\
&=&\gamma_2(L_1,L_2,L)=\langle [L_1,L_2]_2, L\rangle\ .
\end{eqnarray*}
This finishes the proof of (\ref{gum}) and, hence, of the proposition.
\qed

Now we are going to contruct a certain set which should parametrize
the extension equivalence classes for fixed $(\fa,\fl,\rho)$ in the same
way as $H^2(\fl,\fa)$ parametrizes extensions of $\fl$ by the $\fl$-module 
$(\rho,\fa)$.

Since the representation $\rho$ is orthogonal the product 
$\langle\cdot\wedge\cdot\rangle_\fa$
is compatible with the differential
$d$ in the following way
\begin{equation}\label{green}
\langle d\alpha\wedge\tau\rangle_\fa = (-1)^{p+1} \langle \alpha\wedge 
d\tau\rangle_\fa\ ,\qquad \alpha\in C^p(\fl,\fa), \tau \in C^q(\fl,\fa)\ .
\end{equation}
Therefore it induces a kind of cup product on the cohomology groups
\begin{equation}\label{stokes}
\cup: H^p(\fl,\fa)\times H^q(\fl,\fa) \rightarrow 
{\textstyle\bigwedge^{p+q}}\fl^*\ ,\quad [\alpha]\cup [\tau]:=\langle 
\alpha\wedge\tau\rangle_\fa\ .
\end{equation}
Note that the left hand side of (\ref{EK}) is equal to $\frac{1}{2}\langle 
\alpha\wedge\alpha\rangle_\fa$. Thus property
(\ref{EK}) only depends on the cohomology class $a:=[\alpha]\in 
H^2(\fl,\fa)$ 
\footnote{This fact is already implicitly contained in Proposition \ref{abruest}. 
Indeed, its proof shows that if $(\alpha,\gamma)$ defines a metric
Lie algebra, then for {\em any} $\tau\in C^1(\fl,\fa)$ the tupel 
$(a+d\tau,\gamma+\langle(\alpha+\frac{1}{2}d\tau)\wedge\tau\rangle)$ defines a 
metric Lie algebra, too. The Jacobi identity of the latter algebra implies that
$\alpha+d\tau$ satisfies (\ref{EK}).} 
via the condition 
$$
a\cup a =0\ .
$$
Denote by $Z^2_C(\fl,\fa)\subset Z^2(\fl,\fa)$ the space of all cocycles satisfying 
(\ref{EK}).

Since the map $C^1(\fl,\fa)\ni\tau \mapsto \Psi \in GL(\fl^*\oplus\fa\oplus\fl)$ given by (\ref{tsch}) is a group homomorphism up to antisymmetric maps from $\fl$ to $\fl^*$ we can expect to define an action of the abelian group of $1$-cocycles 
$C^1(\fl,\fa)$,
which we write as a right action,
on the direct product $C^2(\fl,\fa)\times {\bigwedge^3}\fl^*$ by
\begin{equation}\label{act}
(\alpha,\gamma)\tau:=(\alpha,\gamma)+\left( d\tau,\langle 
(\alpha+\textstyle\frac{1}{2}d\tau)\wedge\tau \rangle_\fa \right) \ ,\quad \tau \in 
C^1(\fl,\fa),
\alpha\in C^2(\fl,\fa), \gamma \in \textstyle{\bigwedge^3}\fl^*\ .
\end{equation}
Let us verify that (\ref{act}) in fact defines an action. Using (\ref{green})
we observe that for $\tau_1,\tau_2\in C^1(\fl,\fa)$
$$ \langle d\tau_2\wedge \tau_1\rangle_\fa =\langle \tau_2\wedge d\tau_1\rangle_\fa 
=\langle d\tau_1\wedge \tau_2\rangle_\fa\ ,$$
thus 
$$ \textstyle\frac{1}{2}(\langle d\tau_1\wedge\tau_2\rangle_\fa + \langle 
d\tau_2\wedge \tau_1\rangle_\fa)=\langle d\tau_1\wedge\tau_2\rangle_\fa\ .
$$
Therefore we have
\begin{eqnarray*}
(\alpha,\gamma)(\tau_1+\tau_2)\hspace{-3pt}&=&\hspace{-3pt}(\alpha,\gamma)+\left( 
d(\tau_1+\tau_2),\langle 
(\alpha+\textstyle\frac{1}{2}d(\tau_1+\tau_2))\wedge(\tau_1+\tau_2) 
\rangle_\fa\right)\\
\hspace{-3pt}&=&\hspace{-4pt}\left(\alpha+d\tau_1,\gamma+\langle 
(\alpha+\textstyle\frac{1}{2}d\tau_1)\wedge\tau_1)\rangle_\fa\right)
\hspace{-2.4pt}+\hspace{-2.4pt}
\left(d\tau_2,\langle 
(\alpha+d\tau_1+\textstyle\frac{1}{2}d\tau_2)\wedge\tau_2)\rangle_\fa\right)\\
\hspace{-3pt}&=&\hspace{-3pt} ((\alpha,\gamma)\tau_1)\tau_2\ .
\end{eqnarray*}
The subspaces $Z^2_C(\fl,\fa)\times {\bigwedge^3}\fl^*\subset Z^2(\fl,\fa)\times 
{\bigwedge^3}\fl^*\subset C^2(\fl,\fa)\times {\bigwedge^3}\fl^*$ are invariant 
under the action of
$C^1(\fl,\fa)$. The desired parameter
set will be the orbit space
$$ \cH^2_C(\fl,\fa):=(Z^2_C(\fl,\fa)\times 
\textstyle{\bigwedge^3}\fl^*)/C^1(\fl,\fa)\ .$$
We remark that the group of coboundaries $B^1(\fl,\fa)\subset C^1(\fl,\fa)$
acts trivially on $Z^2(\fl,\fa)\times {\bigwedge^3}\fl^*$ because of 
(\ref{green}).
If $(\alpha,\gamma)\in Z^2_C(\fl,\fa)\times {\bigwedge^3}\fl^*$, then we denote 
by $[\alpha,\gamma]$ the corresponding element in 
$\cH^2_C(\fl,\fa)$.

Now the following corollary is just a reformulation of Proposition \ref{abruest}.

\begin{co}\label{triangle}
Fix $(\fa,\fl,\rho)$. The correspondence
$$ \fd_{\alpha,\gamma}(\fa,\fl,\rho)\longmapsto [\alpha,\gamma] \in 
\cH^2_C(\fl,\fa)$$
defines a bijection between the extension equivalence classes of metric
Lie algebras of the form $\fd_{\alpha,\gamma}(\fa,\fl,\rho)$ and elements of
$\cH^2_C(\fl,\fa)$.
\end{co}

\begin{re}{\rm
The set $\cH^2_C(\fl,\fa)$ appears among the cohomology sets introduced and
studied by Grishkov in \cite{grishkov}. There it is denoted by $H^3_\Delta$,
where $\Delta$ is the multiplication given by 
$\frac{1}{2}\langle\cdot\wedge\cdot\rangle_\fa$. This paper also suggests how one 
should define $\cH^2_C(\fl,\fa)$ for non-abelian $\fl$ such that the analogue
of Corollary \ref{triangle} remains true 
(compare Prop. 3.2 of 
\cite{grishkov} which, however, seems to be not quite correct).}
\end{re}

\begin{re}\label{R32}{\rm
Since $\fl$ is abelian the operators $\rho(L)$, $L\in\fl$, have a joint
Jordan decomposition which gives rise to an $\fl$-invariant orthogonal 
decomposition
$$ \fa=\fa^{(\fl)}\oplus \fa^R\ ,$$
where $\fa^{(\fl)}$ is the nil-subspace of the $\fl$-action on $\fa$ and 
$\rho(L)_{|\fa^R}$ is invertible for at least one $L\in \fl$ (in fact for $L$ in an 
open dense subset). As it is well-known, this implies that for all $p$
\begin{equation}\label{coh}
H^p(\fl, \fa^R)=\{0\}\quad  \mbox{and therefore} \quad H^p(\fl, \fa^{(\fl)})\cong 
H^p(\fl, \fa)\ .
\end{equation}
We claim that the embedding 
$$Z^2_C(\fl, \fa^{(\fl)})\times \textstyle{\bigwedge^3}\fl^*
\hookrightarrow Z^2_C(\fl,\fa)\times \textstyle{\bigwedge^3}\fl^* $$
induces an isomorphism
\begin{equation}\label{nl}
\cH^2_C(\fl,\fa^{(\fl)})\cong\cH^2_C(\fl,\fa)\ .
\end{equation}
Indeed, using (\ref{coh}) for $p=1,2$ and that $B^1(\fl,\fa))$ acts trivially   on 
$Z^2(\fl,\fa)\times {\bigwedge^3}\fl^*$ we obtain
\begin{eqnarray*}
\cH^2_C(\fl,\fa)&=& (Z^2_C(\fl,\fa)\times \textstyle{\bigwedge^3}\fl^*)/
C^1(\fl,\fa)\\
&\cong& (Z^2_C(\fl,\fa^{(\fl)})\times \textstyle{\bigwedge^3}\fl^*)/
(C^1(\fl,\fa^{(\fl)})\oplus Z^1(\fl,\fa^R))\\
&=&(Z^2_C(\fl,\fa^{(\fl)})\times \textstyle{\bigwedge^3}\fl^*)/
(C^1(\fl,\fa^{(\fl)})\oplus B^1(\fl,\fa^R))\\
&=&(Z^2_C(\fl,\fa^{(\fl)})\times \textstyle{\bigwedge^3}\fl^*)/
C^1(\fl,\fa^{(\fl)})\\
&=& \cH^2_C(\fl,\fa^{(\fl)})\ .
\end{eqnarray*}
In particular, any extension equivalence class has a 
representative $\fd_{\alpha,\gamma}(\fa,\fl,\rho)$ such that
$\alpha\in Z^2_C(\fl,\fa^{(\fl)})$. 
}
\end{re}

Now we can decide whether two metric Lie algebras 
$\fd_{\alpha_i,\gamma_i}(\fa_i,\fl_i,\rho_i)$, $i=1,2$, are isomorphic
provided that they are regular (see Definition \ref{reg}).
If $S:\fl_1\rightarrow \fl_2$ is a linear map and $\alpha\in C^p(\fl_2,\fa_2)$
then we can form the pull back $S^*\alpha\in C^p(\fl_1,\fa_2)$ given by
$$ S^*\alpha(L_1,\dots,L_p)=\alpha(S(L_1),\dots,S(L_p))\ .$$
We also have the pull back $\bigwedge^p \fl_2^*\ni\gamma\mapsto S^*\gamma\in 
\bigwedge^p \fl_1^*$.
Moreover, cochains can be composed with linear maps $T:\fa_2\rightarrow \fa_1$.
We denote such a composition by $T\alpha$. If $T$ is an isometry
such that $\rho_1(L)=T\rho_2(S(L))T^{-1}$ for all $L\in\fl_1$ and $\alpha\in 
Z^2_C(\fl_2,\fa_2)$,
then $T\,S^*\alpha\in Z^2_C(\fl_1,\fa_1)$. 

\begin{lm}\label{isom1}
We consider the metric Lie algebras
$\fd_{\alpha_i,\gamma_i}(\fa_i,\fl_i,\rho_i)$, $i=1,2$.
If
there exist an isomorphim
$ S:\fl_1\rightarrow \fl_2$
and an isometry
$ U: \fa_1\rightarrow \fa_2 $
such that
\begin{itemize}
\item[(i)] $U\circ\rho_1(L)\circ U^{-1}=\rho_2(S(L)), \ L\in\fl_1,$ and 
\item[(ii)]
$\fd_{\alpha_1,\gamma_1}(\fa_1,\fl_1,\rho_1)$
and 
$\fd_{U^{-1}S^*\alpha_2,S^*\gamma_2}(\fa_1,\fl_1,\rho_1)$ are extension
equivalent, 
\end{itemize}
then $\fd_{\alpha_1,\gamma_1}(\fa_1,\fl_1,\rho_1)$ and 
$\fd_{\alpha_2,\gamma_2}(\fa_2,\fl_2,\rho_2)$ are isomorphic.

Conversely, if $\fd_{\alpha_i,\gamma_i}(\fa_i,\fl_i,\rho_i)$, $i=1,2$, are regular and isomorphic, 
then there exist an isomorphim
$ S:\fl_1\rightarrow \fl_2$
and an isometry
$ U: \fa_1\rightarrow \fa_2 $ satisfying Conditions {\it (i)} and {\it (ii)}.
\end{lm}
\proof
A straightforward computation shows that if $S$ and $U$ satisfy {\it (i)}, then
$S\oplus U\oplus (S^{-1})^*$ defines an isorphism between the metric Lie algebras 
$\fd_{U^{-1}S^*\alpha_2,S^*\gamma_2}(\fa_1,\fl_1,\rho_1)$ and
$\fd_{\alpha_2,\gamma_2}(\fa_2,\fl_2,\rho_2)$.
Thus, if $\Psi: \fd_{\alpha_1,\gamma_1}(\fa_1,\fl_1,\rho_1)\rightarrow
\fd_{U^{-1}S^*\alpha_2,S^*\gamma_2}(\fa_1,\fl_1,\rho_1)$ is an extension 
equivalence, then 
$F:= (S\oplus U\oplus (S^{-1})^*)\circ \Psi$ is an isometric isomorphism between 
$\fd_{\alpha_1,\gamma_1}(\fa_1,\fl_1,\rho_1)$ and 
$\fd_{\alpha_2,\gamma_2}(\fa_2,\fl_2,\rho_2)$.  

In order to prove the opposite direction we really need the regularity assumption. 
We abbreviate $\fd_{\alpha_i,\gamma_i}(\fa_i,\fl_i,\rho_i)$
by $\fd_i$.
Let $F:  \fd_1\rightarrow \fd_2$ be an isomorphism. Then $F(\fz(\fd_1) =\fz(\fd_2)$ 
and $F(\fd_1^\prime)=F(\fd_2^\prime)$. Thus by regularity $F(\fl_1^*)
=\fl_2^*$ and $F(\fa_1\oplus \fl_1^*)=\fa_2\oplus \fl_2^*$. Identifying
the quotients $\fd_i/(\fa_i\oplus \fl_i^*)$ and $(\fa_i\oplus \fl_i^*)/\fl_i^*$
with $\fl_i$ and $\fa_i$, respectively, $F$ induces a bijection $S:\fl_1\rightarrow 
\fl_2$
and an isometry
$ U: \fa_1\rightarrow \fa_2$. Note that $F_{|\fl_1^*}=(S^{-1})^*$.
We compute for $L\in\fl_1, A\in\fa_2$
\begin{eqnarray*}
U\circ\rho_1(L)\circ U^{-1}(A)&=&F([L,U^{-1}(A)])\, \mod \fl_2^*\\
&=& [F(L), F(U^{-1}(A))]\, \mod \fl_2^* = [S(L),A]\, \mod \fl_2^*\\
&=& \rho_2(S(L))\ .
\end{eqnarray*}
Thus $S$ and $U$ satisfy {\it (i)} and, as remarked at the beginning of the proof, 
$S\oplus U\oplus (S^{-1})^*$ is an isomorphism between 
$\fd_{U^{-1}S^*\alpha_2,S^*\gamma_2}(\fa_1,\fl_1,\rho_1)$ and $\fd_2$.
Consequently, $\Psi:=(S^{-1}\oplus U^{-1}\oplus S^*)\circ F$ is an isomorphism
between $\fd_1$ and $\fd_{U^{-1}S^*\alpha_2,S^*\gamma_2}(\fa_1,\fl_1,\rho_1)$
which, by construction, acts on each of the subquotients $\fl_1$, $\fa_1$, and 
$\fl_1^*$ as the identity. Therefore $\Psi$ is an extension equivalence.
This finishes the proof of the lemma.
\qed

Combining Proposition \ref{abruest} with the Lemma \ref{isom1} and identifying 
$(\fa_1,\fl_1)$ with $(\fa_2,\fl_2)$ we obtain

\begin{co}\label{ae}
We consider the metric Lie algebras
$\fd_{\alpha_i,\gamma_i}(\fa,\fl,\rho_i)$, $i=1,2$.
If there exist linear maps $S\in\GL(\fl)$ and
    $U\in O(\fa,\ip_{\fa})$ satisfying
   \begin{equation}\label{ae1}
	U\circ\rho_1(L)\circ U^{-1}=\rho_2(S(L))\ ,  \qquad L\in\fl\ ,
    \end{equation}
    and a linear map $\tau:\fl\longrightarrow \fa$, i.e. $\tau\in
    C^{1}(\fl,\fa)$,
    such that
    \begin{eqnarray}\label{E2}
\hspace{-4em}&&U^{-1}\circ\alpha_{2}(S(L_{1}), S(L_{2}))-\alpha_{1}(L_{1},L_{2})=L_{
1}\tau(L_{2})-L_{2}\tau(L_{1})\\[1ex]
\hspace{-4em}&&\gamma_{2}(SL_{1},SL_{2},SL_{3}) -\gamma_{1}(L_{1},L_{2},L_{3})=
\nonumber\\
\hspace{-4em}&&\qquad\langle\alpha_{1}(L_{1},L_{2}),\tau(L_{3})\rangle_{\fa} +
	\langle\alpha_{1}(L_{3},L_{1}),\tau(L_{2})\rangle_{\fa} +
	\langle\alpha_{1}(L_{2},L_{3}),\tau(L_{1})\rangle_{\fa}\nonumber\\
	\hspace{-4em}&&\qquad +
	\langle \tau(L_{1}),L_{2}\tau(L_{3})\rangle_{\fa}+\langle 
\tau(L_{3}),L_{1}\tau(L_{2})\rangle_{\fa}+
	\langle \tau(L_{2}),L_{3}\tau(L_{1})\rangle_{\fa},\label{E3}
    \end{eqnarray}
    then
    $\fd_{\alpha_{1},\gamma_{1}}(\fa,\fl,\rho_1)$ and
    $\fd_{\alpha_{2},\gamma_{2}}(\fa,\fl,\rho_2)$ are isomorphic.  
In these formulas the action of $\fl$ on $\fa$ is given by $\rho_1$.

Conversely, if $\fd_{\alpha_i,\gamma_i}(\fa,\fl,\rho_i)$, $i=1,2$, are regular and   isomorphic, 
then there exist linear maps $S\in\GL(\fl)$,
    $U\in O(\fa,\ip_{\fa})$, and $\tau\in
    C^{1}(\fl,\fa)$ satisfying Equations {\rm(\ref{ae1}), (\ref{E2})}, and {\rm(\ref{E3})}.
\end{co}

\section{A decomposability criterion}
\label{S4}

There is a natural notion of a direct sum of twofold extensions. 
Namely, $\fd_{\alpha,\gamma}(\fa,\fl,\rho)$ is a direct
sum if and only if there are decompositions $\fl=\fl_1\oplus \fl_2$,
$\fa=\fa_1\oplus\fa_2$, $\fa_1\perp\fa_2$, such that
$\fd_{\alpha,\gamma}(\fa,\fl,\rho)=
\fd_{\alpha_1,\gamma_1}(\fa_1,\fl_1,\rho_1)\oplus
\fd_{\alpha_2,\gamma_2}(\fa_2,\fl_2,\rho_2)$ for appropriate data
$(\rho_i,\alpha_i,\gamma_i)$, $i=1,2$. A decomposable metric Lie algebra
$\fd_{\alpha,\gamma}(\fa,\fl,\rho)$ need not to be a non-trivial direct sum of 
twofold
extensions. However, we have

\begin{lm}\label{predec}
If $\fd_{\alpha,\gamma}(\fa,\fl,\rho)$ 
is extension equivalent to
a non-trivial direct sum of twofold extensions, then it is decomposable. 
Conversely, if $\fd_{\alpha,\gamma}(\fa,\fl,\rho)$ is regular and decomposable, then
it 
is extension equivalent to
a non-trivial direct sum of twofold extensions. 
\end{lm}

\proof
If $\fd_{\alpha,\gamma}(\fa,\fl,\rho)$  
is extension equivalent to
a non-trivial direct sum of twofold extensions, then it is isomorphic 
to a direct sum of two nondegenerate ideals, hence decomposable.

Again, for the opposite direction we need the assumption of regularity.
We write $\fd$ for  $\fd_{\alpha,\gamma}(\fa,\fl,\rho)$ and assume that
there is a non-trivial orthogonal decomposition into ideals 
$$ \fd=\fd_1\oplus \fd_2\ .$$
This induces decompositions
$ \fd^\prime=\fd_1^\prime\oplus \fd_2^\prime$ and 
$\fz(\fd)=\fz(\fd_1)\oplus\fz(\fd_2)$. By regularity we obtain 
$$\fl^*=\fz(\fd)=\fl^*_1\oplus \fl^*_2\ ,$$
where $\fl_i^*:=\fl^*\cap \fd_i=\fz(\fd_i)$, $i=1,2$, and
$$ \fa\oplus\fl^*=\fd_1^\prime\oplus \fd_2^\prime\ .$$
Identifying $\fa$ with the quotient $(\fa\oplus\fl^*)/\fl^*$
we obtain an orthogonal decomposition
$$ \fa=\fa_1\oplus \fa_2\ ,$$
where $\fa_i:=\fd_i^\prime/\fl_i^*$. Note that $\fa_i$, considered
as a subspace of $\fd$, need not to be contained in $\fd_i$.
Similarly
$$ \fl=\fd/\fg^\prime=\fl_1\oplus\fl_2\ ,$$
where $\fl_i:=\fd_i/\fd_i^\prime$. Then the inner product of $\fd_i$ identifies 
$\fl_i^*$ as defined above with the dual of $\fl_i$.
Again, $\fl_i$ considered as a subspace of $\fd$ need not to be contained
in $\fd_i$. However, we have exact sequences
$$ 0\longrightarrow \fa_i \longrightarrow \fd_i/\fl_i^*\longrightarrow 
\fl_i\longrightarrow 0 $$
and
$$ 0\rightarrow \fl_i^* \longrightarrow \fd_i \longrightarrow 
\fd_i/\fl_i^*\rightarrow 0\ ,$$
which give $\fd_i$ the structure of a twofold extension. Hence, after
the choice of appropriate splits, there are isomorphisms 
$$F_i: \fd_i\rightarrow \fd_{\alpha_i,\gamma_i}(\fa_i,\fl_i,\rho_i)$$ 
for certain data
$(\rho_i,\alpha_i,\gamma_i)$ (compare the proof of Proposition \ref{PA}). 
Let $P_i: \fd\rightarrow
\fd_i$ be the projections. Then by construction
$$ (F_1\circ P_1)\oplus (F_2\circ P_2): \fd\rightarrow 
\fd_{\alpha_1,\gamma_1}(\fa_1,\fl_1,\rho_1)\oplus 
\fd_{\alpha_2,\gamma_2}(\fa_2,\fl_2,\rho_2)$$
is an extension equivalence. This finishes the proof of the lemma.
\qed

\begin{pr}\label{Qdec}
    We consider the metric Lie algebra $\fd_{\alpha,\gamma}(\fa,\fl,\rho)$.
    If
    \begin{itemize}
	\item[(i)] there are decompositions
	$$\fa=\fa_{1}\oplus \fa_{2},\quad \fl=\fl_{1}\oplus\fl_{2}$$
	satisfying $\fa_{1}\perp\fa_{2}$ and $\fl_i\oplus \fa_i \ne {0}$, $i=1,2$, 
	and representations $\rho_{i}$ of $\fl_{i}$ on $\fa_{i}$
	such that  $$\rho=\rho_{1}\oplus\rho_{2},$$
	\item[(ii)] the cocycle $\alpha \in Z_C^{2}(\fl,\fa)$ satisfies
	$$\alpha(\fl_{1}\times\fl_{1})\subset\fa_{1},\quad
	\alpha(\fl_{2}\times\fl_{2})\subset\fa_{2},
        $$
	\item[(iii)]  there are linear maps $T_{1}:\fl_{1}\rightarrow\fa_{2}$ 
        and
	$T_{2}:\fl_{2}\rightarrow\fa_{1}$ such that
$$ \alpha(L_1,L_2)=L_1 T_2 L_2-L_2 T_1L_1$$
and
\begin{eqnarray*}
\gamma(L_{1},L_{1}',L_{2})&=&\langle \alpha (L_1,L_1'), T_2(L_2) 
\rangle_\fa +\langle \alpha(L_1',L_2),T_1(L_1)\rangle_\fa\ ,\\
	\gamma(L_{2},L_{2}',L_{1})&=&\langle \alpha (L_2,L_2'), T_1(L_1) 
\rangle_\fa +\langle \alpha(L_2',L_1),T_2(L_2)\rangle_\fa 
\end{eqnarray*}       
	for all $L_{1}, L_{1}'\in\fl_{1}$ and $L_{2}, L_{2}'\in\fl_{2}$
    \end{itemize}
    then $\fd_{\alpha,\gamma}(\fa,\fl,\rho)$ is decomposable.
    Conversely, if $\fd_{\alpha,\gamma}(\fa,\fl,\rho)$ is regular and decomposable, then      
Conditions (i), (ii) and (iii) hold.
\end{pr}
\proof
Let us start the proof with the following remark.
If for $\fd_{\alpha_0,\gamma_0}(\fa,\fl,\rho)$ Conditions {\it (i), (ii)} and {\it (iii)} are 
satisfied with $T_1=T_2=0$, then we will say that $\alpha_0$ and $\gamma_0$ are diagonal. If the 
metric Lie algebra $\fd_{\alpha,\gamma}(\fa,\fl,\rho)$ is extension equivalent to 
$\fd_{\alpha_0,\gamma_0}(\fa,\fl,\rho)$ for diagonal $\alpha_0$ and $\gamma_0$, then it is 
decomposable by Lemma \ref{predec}. Conversely, if $\fd_{\alpha,\gamma}(\fa,\fl,\rho)$ is regular 
and decomposable, then it is extension equivalent to $\fd_{\alpha_0,\gamma_0}(\fa,\fl,\rho)$ with 
diagonal $\alpha_0$ and $\gamma_0$.

Now let us first assume that $\fd_{\alpha,\gamma}(\fa,\fl,\rho)$ is regular and decomposable. Hence 
$\fd_{\alpha,\gamma}(\fa,\fl,\rho)$ is 
extension equivalent to a metric Lie algebra $\fd_{\alpha_0,\gamma_0}(\fa,\fl,\rho)$ with diagonal 
$\alpha_0$ and $\gamma_0$.
By Proposition \ref{abruest} there exists a linear map $\tau:\fl\rightarrow\fa$ 
such
that 
$$ \alpha=\alpha_0+d\tau\quad \mbox{and}\quad \gamma=\gamma_0+\langle (\alpha_0+
\textstyle\frac{1}{2}d\tau)\wedge \tau\rangle_\fa\ .$$
We decompose $\tau=\tau_0+T_1+T_2$, where $\tau_0(\fl_i)\subset \fa_i$,
$T_1(\fl_2)=T_2(\fl_1)=0$, and $T_1(\fl_1)\subset \fa_2$, $T_2(\fl_2)\subset
\fa_1$.
Using that $\rho(L_1)_{|\fa_2}=\rho(L_2)_{|\fa_1}=0$  for $L_i\in\fl_i$ we obtain
$$d\tau(\fl_i\times \fl_i)=d\tau_0(\fl_i\times\fl_i)\subset \fa_i$$
and
$$ d\tau(L_1,L_2)=L_1 T_2(L_2)-L_2 T_1 (L_1)\,,$$
which in particular proves condition {\it (i)}.
It remains to check the last two conditions in {\it (iii)}.
By the above we have
\begin{equation}\label{gugel}
\langle \alpha_0\wedge\tau\rangle_\fa (L_1,L_1',L_2)=\langle \alpha_0 (L_1,L_1'), 
T_2(L_2) \rangle_\fa
\end{equation}
and
\begin{eqnarray}
\textstyle\frac{1}{2}\langle d\tau\wedge \tau\rangle_\fa (L_1,L_1',L_2)
&=&\textstyle\frac{1}{2}\Big(\,\langle d\tau(L_1,L_1'), T_2(L_2)\rangle_\fa
\nonumber\\
&&\quad\label{hupf}
-\langle L_1 T_2(L_2), \tau(L_1') \rangle_\fa + \langle L_2 T_1 (L_1), 
\tau(L_1')\rangle_\fa\\
&&\quad
+\langle L_1' T_2(L_2), \tau(L_1) \rangle_\fa )-\langle L_2 T_1 (L_1'), 
\tau(L_1)\rangle_\fa\Big)\nonumber\\
&=& \textstyle\frac{1}{2}\Big(\,\langle d\tau(L_1,L_1'), T_2(L_2)\rangle_\fa
+\langle T_2(L_2), L_1\tau(L_1')-L_1' \tau(L_1) \rangle_\fa\Big)\nonumber\\
&&- \langle  L_2\tau(L_1'), T_1(L_1) \rangle_\fa\nonumber\\
&=& \langle d\tau(L_1,L_1'), T_2(L_2)\rangle_\fa
+\langle d\tau(L_1',L_2),T_1(L_1)\rangle_\fa\ .\nonumber
\end{eqnarray}
It follows that
\begin{eqnarray*}
\gamma(L_{1},L_{1}',L_{2})&=&\langle (\alpha_0+d\tau) (L_1,L_1'), T_2(L_2) 
\rangle_\fa + \langle d\tau(L_1',L_2),T_1(L_1)\rangle_\fa\\
&=& \langle \alpha (L_1,L_1'), T_2(L_2) 
\rangle_\fa +\langle \alpha(L_1',L_2),T_1(L_1)\rangle_\fa\ .
\end{eqnarray*}
Changing the roles of the indices we also obtain
$$ \gamma(L_{2},L_{2}',L_{1})= \langle \alpha (L_2,L_2'), T_1(L_1) 
\rangle_\fa +\langle \alpha(L_2',L_1),T_2(L_2)\rangle_\fa\ .$$
This finishes the proof of {\it (iii)}.

Vice versa, let us assume that $\fd_{\alpha,\gamma}(\fa,\fl,\rho)$
satisfies the conditions {\it (i)-(iii)} for certain maps $T_1$, $T_2$.
We set 
$$\tau:=T_1+T_2\ ,\quad \alpha_0:= \alpha-d\tau\ ,\quad \mbox{and}\quad
\gamma_0=\gamma-\langle (\alpha_0+
\textstyle\frac{1}{2}d\tau)\wedge \tau\rangle_\fa\ .$$
Then, by Proposition \ref{dreist} the algebras $\fd_{\alpha,\gamma}(\fa,\fl,\rho)$ 
and $\fd_{\alpha_0,\gamma_0}(\fa,\fl,\rho)$
are extension equivalent. By the remark at the beginning of the proof
it is sufficient to show that $\alpha_0$ and $\gamma_0$ are diagonal.
That $\alpha_0$ is diagonal is obvious.
This allows us to repeat the computations (\ref{gugel}) and (\ref{hupf})
in order to obtain
$$ \langle (\alpha_0+
\textstyle\frac{1}{2}d\tau)\wedge \tau\rangle_\fa(L_1,L_1',L_2')=\langle \alpha 
(L_1,L_1'), T_2(L_2) 
\rangle_\fa +\langle \alpha(L_1',L_2),T_1(L_1)\rangle_\fa\ .$$
By assumption, the right hand side is equal to $\gamma(L_1,L_1',L_2)$.
It follows that 
$$\gamma_0(L_1,L_1',L_2)=0\ .$$ 
Changing the roles of the indices, 
we eventually see that $\gamma_0$ is diagonal. This finishes the proof of the 
proposition.
\qed 

We conclude this section with a description of the isomorphism classes
of non-abelian indecomposable metric Lie algebras with the property that $\fg'/\fz(\fg)$ is 
abelian. Let us fix $\fl$, $\fa$, and the inner product $\ip_\fa$.
We consider the set $\Hom(\fl,\so(\fa,\ip_\fa))$ of all orthogonal representations 
of $\fl$ on $\fa$. If $\rho\in \Hom(\fl,\so(\fa,\ip_\fa))$ is fixed
we denote the corresponding $\fl$-module by $\fa_\rho$. 
The group $G:=\GL(\fl)\times O(\fa,\ip_\fa)$ acts on 
$\Hom(\fl,\so(\fa,\ip_\fa))$
by
$(S,U)\rho:=\Ad(U)\circ\rho\circ S^{-1}$, on $C^q(\fl,\fa)$ by
$(S,U)\alpha:=U(S^{-1})^*\alpha$, and on $\bigwedge^3\fl^*$ by $(S,U)\gamma:=(S^{-1})^*\gamma$.
These actions are compatible with the differential $d$, the product 
$\langle\cdot\wedge\cdot\rangle_\fa$, and 
the action (\ref{act}). In particular, if $g\in G$ and $\alpha\in Z^2_C(\fl,\fa_\rho)$, then 
$g\alpha\in Z^2_C(\fl,\fa_{g\rho})$. We
obtain an induced action of $G$ on the disjoint union
$$ \coprod_{\rho\in \Hom(\fl,\so(\fa,\ipa))} \cH^2_C(\fl,\fa_\rho)\ .$$
Let $\cH^2_{C}(\fl,\fa_\rho)_{0}\subset \cH^2_C(\fl,\fa_\rho)$ be the subset
corresponding to all extension equivalence classes of regular indecomposable  
$\fd_{\alpha,\gamma}(\fa,\fl,\rho)$. The relevant conditions for 
$(\alpha,\gamma)$ are given in Lemma \ref{gaga} and Proposition \ref{Qdec}. Then the set
$$ \coprod_{\rho\in \Hom(\fl,\so(\fa,\ipa))} \cH^2_C(\fl,\fa_\rho)_{0}$$ 
is $G$-invariant.
Combining Proposition \ref{PA} and Lemma \ref{isom1} with Corollary \ref{triangle} we obtain
\begin{pr}\label{Pwieder}
    Fix $\fl$, $\fa$, and an inner product $\ip_\fa$. We consider the
    class $\cA(\fl,\fa)$ of non-abelian indecomposable metric Lie algebras $\fg$      
    satisfying
    \begin{enumerate}
        \item $\fg'/\fz(\fg)$ is abelian and isomorphic to $(\fa,\ip_\fa)$ as a    
              semi-Euclidean vector space,
        \item $\dim \fz(\fg)=\dim \fl$.
    \end{enumerate}
    Then the set of isomorphism classes of $\cA(\fl,\fa)$ is in bijective   
    correspondence with the orbit space of the action of $G=\GL(\fl)\times  
    O(\fa,\ip_\fa)$ on 
    $$\coprod_{\rho\in \Hom(\fl,\so(\fa,\ipa))} \cH^2_C(\fl,\fa_\rho)_{0}\ .$$
\end{pr}

\section{The case of maximal isotropic centre}
\label{S5}

In Section \ref{S2} we proved that  any indecomposable metric Lie algebra of signature $(p,q)$, 
$p\le q$, with non-trivial maximal iso\-tropic centre is isomorphic to
 a regular metric Lie algebra $\fd_{\alpha,\gamma}(\fa,\fl,\rho)$ for a Euclidean vector space 
$(\fa,\ip_\fa)$. Therefore we specialize now some of the results from the previous sections to the 
case where $(\fa,\ip_\fa)$ is Euclidean.

\begin{re}\label{Risomax1}
    {\rm
    If $\fa$ is Euclidean the nil-subspace $\fa^{(\fl)}$ coincides 
    with the space $\fa^{\fl}$ of invariants of the representation $\rho$. 
    By Remark \ref{R32} we have 
    $$ {\cal H}^{2}_{C}(\fl,\fa)\cong{\cal H}^{2}_{C}(\fl,\fa^{\fl}).$$
    
    Each extension equivalence class has a representative 
    $\fd_{\alpha,\gamma}(\fa,\fl,\rho)$ such that 
    $\alpha\in Z^{2}_{C}(\fl,\fa^{\fl})$.
    
    Furthermore, we have
    $$
    {\cal H}^{2}_{C}(\fl,\fa^{\fl})=
    \Big( Z^{2}_{C}(\fl,\fa^{\fl})\times\textstyle{\bigwedge^{3}}\fl^{*}\Big)/
    C^{1}(\fl,\fa^{\fl})
    \ = \Big( \textstyle{\bigwedge^{2}_{C}}(\fl^*,\fa^{\fl})\times 
    \textstyle{\bigwedge^{3}}\fl^{*}\Big)/
    C^{1}(\fl,\fa^{\fl})
    $$
    where ${\bigwedge^{2}_{C}}(\fl^*,\fa^{\fl})$ denotes the space of 2-forms on
    $\fl^*$ with values in $\fa^\fl$ satisfying Equation~(\ref{EK}).
    Here the action of $ C^{1}(\fl,\fa^{\fl})$ on 
    $\textstyle{\bigwedge^{2}_{C}}(\fl^*,\fa^{\fl}) \times 
    \textstyle{\bigwedge^{3}}\fl^{*}$
    is given by
    $$(\alpha,\gamma)\tau=(\alpha,\gamma+\langle\alpha\wedge\tau\rangle_{\fa}).$$
    In particular, the orbit of 
    $(\alpha,\gamma)\in\textstyle{\bigwedge^{2}_{C}}(\fl^*,\fa^{\fl}) \times 
    \textstyle{\bigwedge^{3}}\fl^{*}$ is an affine subspace of  
    $\textstyle{\bigwedge^{2}_{C}}(\fl^*,\fa^{\fl}) \times 
    \textstyle{\bigwedge^{3}}\fl^{*}$ with associated vector space 
    $\langle\alpha\wedge C^1(\fl,\fa^\fl)\rangle_\fa \subset 
    {\bigwedge^{3}}\fl^{*}$.
    Hence
    \begin{equation}\label{EN}{\cal H}^{2}_{C}(\fl,\fa^{\fl}) \cong 
    \coprod_{\alpha\in\bigwedge^{2}_{C}(\fl^*,\fa^{\fl})}
    \textstyle{\bigwedge^{3}}\fl^{*}/
    \langle\alpha\wedge C^1(\fl,\fa^\fl)\rangle_\fa \,.
    \end{equation}
    Since ${\cal H}^{2}_{C}(\fl,\fa^{\fl})$ only depends on 
    $\dim\fl$ and $\dim\fa^{\fl}$ we will denote it by 
    ${\cal H}^{2}_{C}(\dim\fl,\dim\fa^{\fl})$.
    }
\end{re}

\begin{pr}\label{Pdeceuc}
    Assume that $(\fa,\ip_{\fa})$ is Euclidean and $\alpha\in 
    \bigwedge^{2}_{C}(\fl^*,\fa^{\fl})$. We also assume that $\fd_{\alpha,\gamma}(\fa,\fl,\rho)$ is 
not two-dimensional (the only two-dimensional algebra of this form is abelian and decomposable).\\
    Then $\fd_{\alpha,\gamma}(\fa,\fl,\rho)$ is decomposable if and 
    only if
    \begin{itemize}
	\item[(i)] there are decompositions
	$$\fa=\fa_{1}\oplus \fa_{2},\quad \fl=\fl_{1}\oplus\fl_{2}$$
	satisfying $\fa_{1}\perp\fa_{2}$ and $\fl_i\oplus \fa_i \ne {0}$, $i=1,2$, 
	and representations $\rho_{i}$ of $\fl_{i}$ on $\fa_{i}$
	such that  $$\rho=\rho_{1}\oplus\rho_{2},$$
	\item[(ii)] the cocycle $\alpha \in \bigwedge^{2}_{C}(\fl^*,\fa^{\fl})$ satisfies
	$$\alpha(\fl_{1}\times\fl_{1})\subset\fa_{1},\quad 
	\alpha(\fl_{2}\times\fl_{2})\subset\fa_{2},
        \quad \alpha(\fl_{1}\times\fl_{2})=0,$$
	\item[(iii)]  there are linear maps $T_{1}:\fl_{1}\rightarrow\fa_{2}$ and
	$T_{2}:\fl_{2}\rightarrow\fa_{1}$ such that
	$$\gamma(L_{1},L_{1}',L_{2})=\langle\alpha(L_{1},L_{1}'), 
	T_{2}L_{2}\rangle_{\fa},\quad 
	\gamma(L_{2},L_{2}',L_{1})=\langle\alpha(L_{2},L_{2}'), 
	T_{1}L_{1}\rangle_{\fa}$$
	for all $L_{1}, L_{1}'\in\fl_{1}$ and $L_{2}, L_{2}'\in\fl_{2}$.
    \end{itemize}
\end{pr} 
\proof If Conditions {\it (i), (ii), (iii)} are satisfied, then $\fd_{\alpha,\gamma}(\fa,\fl,\rho)$ 
is decomposable by Proposition \ref{Qdec}. 

Assume now that $\fd_{\alpha,\gamma}(\fa,\fl,\rho)$ is decomposable. If 
$\fd_{\alpha,\gamma}(\fa,\fl,\rho)$ is regular, then we can again apply Proposition \ref{Qdec} and 
the assertion follows since 
$$\alpha(L_1,L_2)=L_1 T_2L_2-L_2T_1L_1\in\fa^\fl$$ implies $\alpha(L_1,L_2)=0$
for $L_1\in\fl_1$ and $L_2\in\fl_2$. If $\fd_{\alpha,\gamma}(\fa,\fl,\rho)$ is not regular then 
there exist elements $A_0\in\fa$, $L_0\in\fl$, $A_0+L_0\not=0$ which satisfy Equations 
(\ref{sing}). If $L_0=0$, then 
$$\fa_1=\RR A_0,\ \fa_2=(\RR A_0)^\perp,\ \fl_1=0,\ \fl_2=\fl,\ T_1=T_2=0$$
satisfy Conditions {\it (i), (ii), (iii)}. For $L_0\not=0$ we can take  
$$\fa_1=0,\ \fa_2=\fa,\ \fl_1=\RR L_0,\  T_1(L_0)=A_0,\ T_2=0$$
and choose for $\fl_2$ any complement of $\RR L_0$ in $\fl$. This decomposition
of $\fl\oplus\fa$ is non-trivial in the sense of Condition {\it (i)} since
by assumption $\fd_{\alpha,\gamma}(\fa,\fl,\rho)$ is at least 
three-dimensional.  
\qed
\begin{re}\label{Risomax2}
    {\rm
    Let $\fl$ be an abelian Lie algebra of dimension $l$ and 
    $(\fa,\ip_{\fa})$ a Euclidean space of dimension $n$. Recall that 
    $\GL(\fl)\times O(\fa,\ip_{\fa})$ acts on 
    $\Hom(\fl,\so(\fa,\ip_{\fa})$ by
    $$(S,U)\rho (\cdot) = U\circ \rho (S^{-1}(\cdot))\circ U^{-1}.$$
    The orbit space of this action equals
    \begin{eqnarray*} 
	\Hom(\fl,\so(\fa,\ip_{\fa})/ \GL(\fl)\times O(\fa,\ip_{\fa})&=& 
        (\fl^*)^{m}/(\GL(\fl)\times\frak 
        S_{m}\ltimes(\ZZ_{2})^{m}) \\
	&= &\coprod_{k\le l}{\Bbb G}_{k,m}/(\frak 
        S_{m}\ltimes(\ZZ_{2})^{m}),
    \end{eqnarray*}
    where
    $m=\left[\frac n2\right]$, ${\Bbb G}_{k,m}$ 
    denotes the Grassmannian and $\frak S_{m}$ the symmetric group. Here we first 
    identify the orbit of $\rho\in\Hom(\fl,\so(\fa,\ip_{\fa})$ with the 
    $\GL(\fl)\times\frak S_{m}\ltimes(\ZZ_{2})^{m}$-orbit of the weights 
    $\lambda=(\lambda^1,\dots,\lambda^m)\in(\fl^*)^m$ of $\rho$ and then with the 
    $\frak S_{m}\ltimes(\ZZ_{2})^{m}$-orbit of $\Span\{\lambda(L)\mid L\in\fl\}$.
    }
\end{re}
\begin{ex}
    {\rm
    Consider the abelian Lie algebra 
    $\fl=\RR^{l}$ 
    and let $(\fa,\ip_{\fa})$ be the Euclidean space $\RR^{2m}\oplus\RR^{k}$.
    Choose 
    $\lambda=(\lambda^1,\dots,\lambda^m)\in(l^*)^m$.
    Define an orthogonal representation $\rho$ of $\fl$ on $\fa$ by
    $$\rho(L)=\left(
    \begin{array}{ccc}
        0 & -\diag(\lambda(L)) &0 \\
        \diag(\lambda(L)) & 0&0\\
        0&0&0
    \end{array} \right).$$
    Furthermore, assume $\alpha\in \bigwedge^{2}_{C}(\RR^{l},\RR^{k})$, 
    $\gamma\in\bigwedge^{3}\RR^{l}$ and define
    $$\fd_{\alpha,\gamma}(m,k,l,\lambda) 
    :=\fd_{\alpha,\gamma}(\RR^{2m+k},\RR^{l},\rho).$$
    }
\end{ex}    
Next we want to give a method which yields (in principle) a description of the 
isomorphism classes of indecomposable metric Lie algebras of 
signature $(l,n+l)$, $l>0$, with maximal isotropic centre. By Theorem \ref{T} this set is equal 
to the set of isomorphism classes in ${\cal A}(\fl,\fa)$ for a fixed $n$-dimensional 
Euclidean space $(\fa,\ip_{\fa})$ and a fixed $l$-dimensional abelian 
Lie algebra $\fl$. By Proposition \ref{Pwieder} the isomorphism classes 
in ${\cal A}(\fl,\fa)$ are in bijective 
correspondence with the $\GL(\fl)\times O(\fa,\ip_{\fa})$--orbits of
$$\coprod _{\rho\in\Hom(\fl,\so(\fa,\ipa))}{\cal 
H}^{2}_{C}(\fl,\fa_{\rho})_{0}\,.$$ Note that in the case of Euclidean $\fa$ the set
${\cal H}^{2}_{C}(\fl,\fa_{\rho})_{0}$ consists of all elements in ${\cal 
H}^{2}_{C}(\fl,\fa_{\rho})$ which correspond to indecomposable metric Lie algebras.
According to Remark \ref{Risomax1} and Remark \ref{Risomax2} this 
orbit space equals
\begin{eqnarray*}
    && \coprod_{k=0}^n\Big(\Big(\coprod _
    {\rho\in\Hom(\fl,\so(\frak a,\ipa))\atop
    \dim \fa^\fl_\rho=k}{\cal 
    H}^{2}_{C}(\fl,\fa_{\rho})_{0}\Big) 
    \Big/O(\fa,\ip_{\fa})\Big)\Big/ \GL (\fl)\\
    &=&\coprod_{\scriptsize{\begin{array}{c}(k,m)\\ k+2m=n\end{array}}}\left(\ 
\Big((\fl^{*}\setminus 0)^{m}
    /\frak 
    S_{m}\ltimes(\ZZ_{2})^{m}\Big)\times\Big({\cal 
    H}^{2}_{C}(l,k)/O(k)\Big) \ \right)_{0}\Big/ \GL(\fl),
\end{eqnarray*}	
where the subscript $0$ again denotes the subset of elements which correspond to indecomposable 
metric Lie algebras. By Proposition \ref{Pdeceuc} indecompasability implies that 
$\alpha(\fl,\fl)=a_\rho^\fl$. In fact, this property is
already a consequence of regularity. In particular, ${\cal H}^{2}_{C}(\fl,\fa_{\rho})_{0}=0$ for 
$k=\dim a_\rho^\fl>l(l-1)2$, $l=\dim \fl$.\\

Hence, if we want to determine the isomorphism classes in ${\cal A}(\fl,\fa)$, then for all $(k,m)$ 
with $k+2m=n$ we have to
\begin{enumerate}
    \item determine ${\cal H}^{2}_{C}(l,k)$ for $k\le l(l-1)/2$ according to (\ref{EN}),
    \item determine representatives 
          $[\alpha_{\iota},\gamma_{\iota}]$ ($\iota$ runs over an index set $I$) 
          of those 
          $(\GL(\fl)\times O(k))$-orbits in ${\cal H}^{2}_{C}(l,k)$ 
          whose elements $[\alpha,\gamma]$ satisfy 
          $\alpha(\fl,\fl)=\RR^{k}$, 	  
    \item determine the stabilizer $G_{\iota}\subset \GL(\fl)\times O(k)$ 
          for each $[\alpha_{\iota},\gamma_{\iota}]$,
    \item determine the subset $\Lambda\subset (\fl^{*}\setminus 0)^{m}$
          of those $\lambda\in(\fl^{*}\setminus 0)^{m}$ for which 
          $\fd_{\alpha_{\iota},\gamma_{\iota}}(m,k,l,\lambda)$ is indecomposable 
          according to Proposition \ref{Pdeceuc},  
    \item determine the set $\{{\cal O}_\theta\mid\theta\in\Theta_\iota\}$, of $\bar    
          G_{\iota}\times \frak S_{m}\ltimes(\ZZ_{2})^{m}$-orbits in $\Lambda$, where
	  $\bar G_{\iota}={\rm proj}_{\GL(\fl)}G_{\iota} \subset \GL(\fl).$
\end{enumerate} 

Then each isomorphism class of indecomposable metric Lie algebras of 
signature $(l,n+l)$ with maximal isotropic centre is represented by 
one of the metric Lie algebras 
$\fd_{\alpha_{\iota},\gamma_{\iota}}(m,k,l,\lambda)$,
where $2m+k=n$. Here $\iota\in I$ and the $\bar G_{\iota}$-orbit  
${\cal O}_\theta$, $\theta\in\Theta_{\iota}$ of $\lambda$ are uniquely determined.\\

\begin{re}\label{dick}
{\rm
The main obstruction (besides problems of bookkeeping and the lack of time)
which prevents us from carrying out this procedure
in general is the rather complicated geometry of 3-forms. Indeed, Step 2 involves
the determination of the orbit spaces
$$  \left(\textstyle\bigwedge^{3}\fl^{*}/
    \langle\alpha\wedge C^1(\fl,\RR^{k})\rangle\right)/\bar G_\alpha\ ,\quad \alpha\in 
\textstyle\bigwedge^{2}_{C}(\fl^*,\RR^{k})
   \mbox{ such that }\alpha(\fl,\fl)=\RR^{k}\ ,$$
where $\bar G_\alpha\subset GL(\fl)$ denotes projection of the stabilizer of $\alpha$ in
$GL(\fl)\times O(k)$ onto the first factor.
This problem is already difficult in the case $\alpha=0$ (which implies $k=0$).
To our knowledge, the relevant orbit space
$$ \textstyle\bigwedge^{3}\fl^{*}/GL(\fl) $$ 
is not known for $l\ge 10$. If $l\le 8$, then it is finite; see \cite{gurevich}, \S 35 for a 
classification over $\CC$, compare also \cite{hitchin00}.
For dimension $l=9$ we refer to \cite{katanova} and the references cited therein. Here the orbit 
space is (generically) 3-dimensional.}
\end{re}

Now we will apply this method for metrics with small index. Let us start with Lorentzian Lie 
algebras. We already know from Corollary \ref{i2} that any indecomposable non-simple Lorentzian Lie 
algebra has a maximal isotropic centre. Therefore our recipe reproduces the known classification of 
(at least 2-dimensional) indecomposable non-simple Lorentzian Lie algebras (see \cite{M85}). Any 
such Lie algebra is isomorphic to exactly one of the indecomposable Lorentzian algebras
$$\osc(\lambda^1,\dots,\lambda^m):=\fd_{0,0}(m,0,1,\lambda),\quad 
1=\lambda^1\le\lambda^2\le\dots\le\lambda^m,$$ 
where $(\lambda^1,\dots,\lambda^m)=\lambda(L)$ for a fixed basis $L$ of $\fl=\RR$. 

Note, that the only simple Lorentzian Lie algebra is $\fsl(2,\RR)$ since the Killing forms of all 
other simple Lie algebras (up to the sign) do not have Lorentzian signature. This completes the 
classification of all indecomposable Lorentzian Lie algebras.

In the following we will fix a basis $L_1,\dots,L_l$ of $\fl=\RR^l$, the dual basis $Z_1,\dots,Z_l$ 
in $\fl^*$, and an orthonormal basis $A_1,\dots,A_k$ of $\RR^k$. For $\lambda\in (\fl^*)^m$ we 
define $(\lambda^i_j)^{i=1,\dots,m}_{j=1,\dots,l}=(\lambda^i(L_j))^{i=1,\dots,m}_{j=1,\dots,l}$ and 
denote by
$\lambda^i=(\lambda^i_1,\dots,\lambda^i_l)$, $i=1,\dots,m$, and 
$\lambda_j=(\lambda^1_j,\dots,\lambda^m_j)^\top$, $j=1,\dots,l$ the rows and the columns of 
$(\lambda^i_j)^{i=1,\dots,m}_{j=1,\dots,l}$, respectively.

\begin{theo}\label{iso23}
For $l=2,3$ each isomorphism class of indecomposable metric Lie algebras of 
signature $(l,l+n)$ with maximal isotropic centre is represented by 
a metric Lie algebra $\fd_{\alpha_{\iota},\gamma_{\iota}}(m,k,l,\lambda)$, $\lambda\in 
\Lambda_\iota$,
for exactly one combination of $m,\,k$, $\alpha_\iota$, $\gamma_\iota$ in the table below and for a 
suitable $\lambda$ in the set $\Lambda_\iota\subset (\fl^*\setminus 0)^m$, which is also given in 
the table. 
In particular, we have $n=2m+k$.

The metric Lie algebras $\fd_{\alpha_{\iota},\gamma_{\iota}}(m,k,l,\lambda)$ and 
$\fd_{\alpha_{\iota},\gamma_{\iota}}(m,k,l,\bar\lambda)$ are isomorphic if and only~if
$$\vartheta_\iota(\lambda)=\vartheta_\iota(\bar\lambda)\ \mod \frak S_{m}\ltimes(\ZZ_{2})^{m}$$
for the map $\vartheta_\iota$ which is associated with the combination $m,\,k$, $\alpha_\iota$, 
$\gamma_\iota$ and listed in the table.

Conversely, all the metric Lie algebras $\fd_{\alpha_{\iota},\gamma_{\iota}}(m,k,l,\lambda)$ with 
$l,m,k$,
$\alpha_{\iota},\gamma_{\iota}$, $\lambda\in\Lambda_\iota$ as in the table are indecomposable and 
of 
signature $(l,l+2m+k)$ with maximal isotropic centre.\\

\hspace{0.1em} 
\fbox{ $l=2$ }\\[0.5ex]
{\rm
\begin{center}
\begin{tabular}{|l|l|}\hline &\\[-2ex]
$k=0$ & \parbox{12.3cm}{$\alpha=0$, $\gamma=0$}\\
$m\ge 3$ &\\[-1.5ex]
      & $\Lambda=\left\{ \lambda\in (\fl^*\setminus 0)^m \left| \begin{array}{l}
	\{\lambda^i\mid i=1,\dots,m\}\subset\fl^*\ \mbox{\rm
        is not contained in the} \\ \mbox{\rm union of two 1-dimensional subspaces of 
        } \fl^*
	\end{array}\right.\right\}$\\&\\[-1.5ex] 
  & $\vartheta(\lambda)=\Span\{\lambda_1,\lambda_2\}\in {\Bbb G}_{2,m}$	        
\\[2ex]
\hline
\end{tabular}
\end{center}
\newpage 
\begin{center}
\begin{tabular}{|l|l|}
\hline &\\[-2ex]
$k=1$ &  \parbox{12.3cm}{$\alpha=Z_{1}\wedge Z_{2}\otimes A_{1}$ , $\gamma=0$}\\
$m\ge 0$ &\\[-1.5ex]
	& $\Lambda=(\fl^*\setminus 0)^m$\\ &\\[-1.5ex]
       & $\vartheta(\lambda)=(E_\lambda,\omega_\lambda),\ $  where 
         $E_\lambda=\Span\{\lambda_1,\lambda_2\}\subset \RR^m$,\\ 
       & \hspace{3.92cm}  $\ \omega_\lambda = [\lambda_1\wedge\lambda_2]\,\in 
         \bigwedge^2(E_\lambda)/\{\pm1\}$
\\[2ex]
\hline
\end{tabular}
\end{center}
\vspace{0.7cm}
\hspace{0.1em} \fbox{ $l=3$ }\\[0.5ex]
\begin{center}
\begin{tabular}{|l|l|}\hline &\\[-2ex]
$k=0$ & \parbox{12.3cm}{$\alpha_1=0$, $\gamma_1=0$}\\
$m\ge 4$ &\\[-1.5ex]
      & $\Lambda_1=\left\{ \lambda\in (\fl^*\setminus 0)^m \left| \begin{array}{l}
	\{\lambda^i\mid i=1,\dots,m\}\subset\fl^*\ \mbox{\rm
        is not contained  } \\ \mbox{\rm  in the union of a 2-dimensional and a }
        \\ \mbox{\rm 1-dimensional subspace of  
        } \fl^*
	\end{array}\right.\right\}$\\&\\[-1.5ex]
      & $\vartheta_1(\lambda)=\Span\{\lambda_1,\lambda_2,\lambda_3\}
	\in {\Bbb G}_{3,m}$
\\[2ex]
\hline&\\[-2ex]
$k=0$ & $\alpha_2=0$, $\gamma_2=Z_{1}\wedge Z_{2}\wedge Z_{3}$\\
$m\ge 0$ & \\[-1.5ex]
      & $\Lambda_2=(\fl^*\setminus 0)^m$\\ &\\[-1.5ex]
      & $\vartheta_2(\lambda)=(E_\lambda,\omega_\lambda),\ $ where $E_\lambda=   
      \Span\{\lambda_1,\lambda_2,\lambda_3\}\subset \RR^m,$\\
      & \hspace{4.21cm} $\omega_\lambda= \lambda_1\wedge\lambda_2\wedge\lambda_3\in 
      \bigwedge^3(E_\lambda)$	
\\[2ex]\hline 
 &\\[-2ex]
$k=1$ & $\alpha=Z_{1}\wedge Z_{2}\otimes A_{1}$ , $\gamma=0$\\
$m\ge2$ &\\[-1.5ex]
      & $\Lambda=\left\{ \lambda\in (\fl^*\setminus 0)^m \left|\:\dim\Span\{\lambda^j\mid 
\lambda^j_3\not=0\}>1\right.\right\}$\\ &\\[-1.5ex]
      & $\vartheta(\lambda)=(E_\lambda,F_\lambda,\omega_\lambda),\ $ where 
        $E_\lambda=\Span\{\lambda_1,\lambda_2,\lambda_3\}\subset \RR^m,$\\ 
      & \hspace{4.78cm}$F_\lambda=\RR \lambda_3\subset E_\lambda,$\\
      & \hspace{4.79cm}$\omega_\lambda= [\lambda_1\wedge\lambda_2 \mod \RR \lambda_3]\in 
        \bigwedge^2(E_\lambda/F_\lambda)/\{\pm1\}$
\\[2ex] 
\hline&\\[-2ex]
$k=2$ & \parbox{12cm}{$\alpha=Z_{1}\wedge Z_{2}\otimes A_{1}+Z_{1}\wedge Z_{3}\otimes A_{2}$, 
        $\gamma=0$}\\
$m\ge0$ &\\[-1.5ex]
	& $\Lambda=(\fl^*\setminus 0)^m$\\ &\\[-1.5ex]
        & $\vartheta(\lambda)=(B_\lambda,v_\lambda)\mod \RR^* ,$\\
        & \hspace{1.8cm}where $ B_\lambda=\Big(\,(\lambda^i_2,\lambda^i_3) 
          \cdot(\lambda^j_2,\lambda^j_3)^\top\Big)^{i=1,\dots,m} 
          _{j=1,\dots,m}\in \End(\RR^m),$\\
        &\hspace{3cm}$v_\lambda=[\lambda_1]\in\Coker B_\lambda,$\\[1ex]
        & \hspace{3cm}$c\in\RR^*$ acts by  
        $(B,v)c=(c^2B,\frac{1}{c}v)$
\\[2ex]
\hline &\\[-2ex]   
$k=3$ & $\alpha=Z_{1}\wedge Z_{2}\otimes A_{1}+Z_{1}\wedge Z_{3}\otimes A_{2}
        +Z_{2}\wedge Z_{3}\otimes A_{3}$ , $\gamma=0$\\
$m\ge0$ &\\[-1.5ex]        
        & $\Lambda=(\fl^*\setminus 0)^m$\\ &\\[-1.8ex]
        & $\vartheta(\lambda)=\Big(\lambda^i  
        \cdot(\lambda^j)^\top\Big)^{i=1,\dots,m} 
        _{j=1,\dots,m}\in \End(\RR^{m})$     
\\[2ex]        
\hline
\end{tabular}
\end{center}
}
\end{theo}
\newpage
\proof
We restrict ourselves to the case $l=3$. We follow the recipe described above. First we note that 
for $l=3$ Equation (\ref{EK}) is trivially satisfied. Then it is not hard to see that each class in 
${\cal H}^{2}_{C}(3,k)$ is represented by exactly one element 
$$(\alpha,\gamma) \in \Big(\{0\}\times\textstyle{\bigwedge}^3\fl^*\Big)\cup 
\Big(\,(\textstyle{\bigwedge} ^2(\fl^*,\RR^k)\setminus\{0\})\times\{0\}\Big)\subset 
\textstyle{\bigwedge} ^2(\fl^*,\RR^k)\oplus \textstyle{\bigwedge}^3\fl^*\,.$$

The group $\GL(\fl)\times O(k)$ acts on this set of representatives. We have to determine the 
orbits of this action and to choose representatives, where we may restrict ourselves to orbits of 
those $(\alpha,\gamma)$ which satisfy $\dim (\alpha(\fl,\fl))=k$. Obviously, for $k=0$, i.e. 
$\alpha=0$ we have two such orbits, one orbit consists of the element $(\alpha,\gamma)=(0,0)$, the 
other one is characterized by $\gamma\not=0$.
For $k\ge 1$ we have $\alpha\not=0$ and may therefore assume that $\gamma=0$ holds. We claim that 
in this case there is only one orbit. The representation $\bigwedge^2\fl$ of $\GL(\fl)$ is 
equivalent to the representation $\rho$ of $\GL(\fl)$ on $\fl^*$ given by 
$\GL(\fl)\ni S \mapsto \det S\cdot (S^{-1})^*\in\GL(\fl^*)$. In particular, $\GL(\fl)$ acts 
transitively on oriented bases in $\bigwedge^2\fl$. Therefore, given a basis 
$\sigma^1,\sigma^2,\sigma^3$ in $\bigwedge^2\fl$, a basis $A_1,\dots,A_k$ of $\RR^k$ and a 2-form 
$\alpha\in\bigwedge^2(\fl^*,\RR^k)$ we find a map $S\in\GL(\fl)$ such that 
$(S^*\alpha)(\sigma^j)=\pm A_j$ for $1\le j\le k$ and $(S^*\alpha)(\sigma^j)=0$ for $k< j\le 3$. 
Now it is obvious that there exists a map $U\in O(k)$ such that $(US^*\alpha)(\sigma^j)= A_j$ for 
$1\le j\le k$ and $(US^*\alpha)(\sigma^j)=0$ for $k< j\le 3$, which proves the claim. In 
particular, we may choose 
$(\alpha,\gamma)$ as in the table.

Since the proofs of the assertions on indecomposability and isomorphy are similar for different $k$ 
we will give here only the one for $k=2$, which is in some sense the most complicated one. First we 
determine the projection $\bar G={\rm proj}_{\GL(\fl)}G$ of the stabilizer $G\subset\GL(\fl)\times 
O(2)$ of $(\alpha,\gamma)=(Z_{1}\wedge Z_{2}\otimes A_{1}+Z_{1}\wedge Z_{3}\otimes A_{2},0)$ on 
$\GL(\fl)$:
\begin{eqnarray}\label{gleich}
\bar G &=&\left\{S\in\GL(\fl)\left| 
      \begin{array}{l} 
          S^*\alpha(L_1,L_2), S^*\alpha(L_1,L_3)\mbox{ is an orthonormal basis of } 
          \RR^k,\\
          S^*\alpha(L_2,L_3)=0
      \end{array}
      \right.
      \right\}\nonumber\\
      &=&\left\{S\in\GL(\fl)\left|
      \begin{array}{l} 
          S(\Span\{L_2,L_3\})=\Span\{L_2,L_3\}\\
          S|_{\Span\{L_2,L_3\}}\in c\cdot O(2) \ \mbox{with respect to  } \\
          \mbox{the basis }\{L_2,L_3\}, \mbox { where } c=| \det S|
      \end{array}
      \right.
      \right\}.
\end{eqnarray}

Next we determine the set $\Lambda\subset (\fl^{*}\setminus 0)^{m}$
of those $\lambda\in(\fl^{*}\setminus 0)^{m}$ for which 
$\fd_{\alpha,0}(m,2,3,\lambda)$ is indecomposable. Assume that $\fd_{\alpha,0}(m,2,3,\lambda)$ is 
decomposable. By Proposition \ref{Pdeceuc} we have decompositions $\fl=\RR^3=\fl_1\oplus\fl_2$ and 
$\fa=\RR^{2m+k}=\fa_1\oplus\fa_2$, at least one of them non-trivial, such that $\rho=\rho_1\oplus 
\rho_2$ holds for the representation $\rho$ defined by $\lambda\in(\fl^{*}\setminus 0)^{m}$ and 
$\alpha$ satisfies $\alpha(\fl_{1}\times\fl_{1})\subset\fa_{1}$, 
$\alpha(\fl_{2}\times\fl_{2})\subset\fa_{2}$, $\alpha(\fl_{1}\times\fl_{2})=0$. Since $\Ker 
\alpha\subset \bigwedge^2(\fl^*,\fa)$ is only 1-dimensional the condition 
$\alpha(\fl_1\times\fl_2)=0$ implies that the decomposition of $\fl$ is trivial.
Hence, we may assume $\fl_1=\fl$ and $\fl_2=0$. Thus the decomposition of $\fa$ is non-trivial, 
i.e. $\fa_2\not=0$. Now $\alpha(\fl_{1}\times\fl_{1})\subset\fa_{1}$ implies 
$\alpha(\fl\times\fl)=\RR^k\subset\fa_1$. On the other hand $\fl_2=0$ yields 
$0\not=\fa_2\subset\fa^\fl$ which is a contradiction to $\lambda^j\not= 0$ for $j=1,\dots,m$. 
Consequently, $\Lambda=(\fl^{*}\setminus 0)^{m}$.

Finally, we determine the $\bar G$-orbits in $(\fl^{*}\setminus 0)^{m}$.  From (\ref{gleich}) we 
know that $(\lambda_1,\lambda_2,\lambda_3)$ and $(\bar\lambda_1,\bar\lambda_2,\bar\lambda_3)$ are 
in the same orbit if and only if there is a map $S'\in O(2)$ and a real number $c\in\RR^*$ such 
that 
$c(\lambda_2,\lambda_3)\cdot S'=(\bar\lambda_2,\bar\lambda_3)$ and $\lambda_1=c\bar\lambda_1\ \mod\ 
\Span\{\lambda_2,\lambda_3\}$. Now we use the following general fact, which is true for arbitrary 
homomorphisms $A,\bar A:\RR^N\rightarrow \RR^m$:
$$\exists s\in O(N): \bar A=A\cdot s \quad \Leftrightarrow \quad A\, A^*=\bar A \, \bar 
A^*:\RR^m\longrightarrow \RR^m. $$
We apply this for $A=(\lambda_j^i)^{i=1,\dots,m}_{j=2,3}$ and $\bar 
A=(\bar\lambda_j^i)^{i=1,\dots,m}_{j=2,3}$ and use the notation 
$B_\lambda=(\lambda_j^i)^{i=1,\dots,m}_{j=2,3} ((\lambda_j^i)^{i=1,\dots,m}_{j=2,3})^\top$. Note 
that $\Span \{\lambda_2,\lambda_3\}=\im (\lambda_j^i)^{i=1,\dots,m}_{j=2,3} = \im B_\lambda$. Now 
we can say that $(\lambda_1,\lambda_2,\lambda_3)$ and $(\bar\lambda_1,\bar\lambda_2,\bar\lambda_3)$ 
are in the same $\bar G$-orbit if and only if there exists a real number $c\in\RR^*$ such that
$ c^2B_\lambda= B_{\bar\lambda}$ and $[\lambda_1]=c[\bar\lambda_1]\in \Coker B_\lambda$ $(=\Coker 
B_{ \bar\lambda}).$
 \qed
 
 \begin{re}
{\rm 
Based on the classification results mentioned in Remark \ref{dick}, it should be not too difficult 
to obtain a classification of indecomposable metric Lie algebras with maximal isotropic
centre of signature $(l,l)$, $2\le l\le 9$. The isomorphism classes of these
metric Lie algebras are parametrized by those $[\gamma] \in \bigwedge^{3}\fl^{*}/GL(\fl)$ which are 
not represented by elements of the form
$\gamma=\gamma_1+\gamma_2$, $\gamma_i\in \bigwedge^{3}\fl_i^{*}$ for a nontrivial
decomposition $\fl=\fl_1\oplus\fl_2$, via
$$ [\gamma]\longmapsto \fd_{0,\gamma}(0,0,l,0)\ .$$  
In fact, it is easy to carry out the classification for $l\le 5$. There are no such algebras
for $l=2,4$ and exactly one for $l=3,5$. Moreover, it can be shown that there are exactly two 
for $l=6$.}
\end{re}


\section{Metric Lie algebras of index 2}

\begin{ex} \label{ExdA}
    {\rm
    For a map $A\in\so(p,q)$ we consider 
    $\fd_A(\RR^{p,q},\RR):=\fd_{0,0}(\RR^{p,q},\RR,\rho)$, where $\rho$ 
    is defined by $\rho(1)=A$. Recall that $\fd_A(\RR^{p,q},\RR)$ is a metric Lie 
    algebra of signature $(p+1,q+1)$ which can be described as follows:
    $\fd_A(\RR^{p,q},\RR)$ is isomorphic to $\RR^{1,1}\oplus\RR^{p,q}$ as a  
    pseudo-Euclidean vector 
    space, there exists a basis 
    $L,Z$ of $\RR^{1,1}$ with $\langle L,L\rangle = \langle Z,Z\rangle=0$ and
    $\langle L,Z\rangle=1$ such that $Z\in\fz(\fd_A(\RR^{p,q},\RR))$ and   
    $$ [X_1,X_2]=\langle AX_1,X_2\rangle Z,\quad [L,X]=A(X) $$
    for all $X\in\RR^{p,q}$. 
    }
\end{ex}    
\begin{ex}{\rm
    Consider the metric Lie algebra 
    $\osc(\lambda_1,\lambda_2) = \fd_{0,0}(m,0,2,\lambda)$ of 
    signature $(2,2+2m)$ which was defined in the previous section, where
    $\lambda_1=\lambda(L_1),\ \lambda_2=\lambda(L_2)\in\RR^m$ for a fixed basis $L_1,L_2$ of 
    $\RR^2$.  Recall that 
    this algebra can be described in the following way: $\osc(\lambda_1,\lambda_2)$ 
    for $\lambda_1,\,\lambda_2\in\RR^m$
    is isomorphic to $\RR^{2,2}\oplus\RR^{2m}$ as a pseudo-Euclidean vector space, 
    there exists a basis $L_1, L_2, Z_1, Z_2$ of $\RR^{2,2}$ with
    $$\langle L_i,L_j\rangle = \langle Z_i,Z_j\rangle=0,\ \langle Z_i,L_j\rangle
    =\delta_{ij}$$ for $i,j=1,2$ and an orthonormal basis 
    $X_{1},\dots,X_{m},Y_{1},\dots,Y_{m}$ of $\RR^{2m}$ such that
    $$\begin{array}{l}
        \, Z_{1},Z_{2}\in\fz(\osc(\lambda_1,\lambda_2))  \\[0.5ex]
        \,[X_{j},Y_{k}]=\delta_{jk}(\lambda_1^{j}Z_{1}+
	\lambda_2^{j}Z_{2}),
	\quad [X_{j},X_{k}]=0,\quad [Y_{j},Y_{k}]=0   \\[0.5ex]
        \,[L_{i},X_{j}]=\lambda_i^{j}Y_{j},\quad
        \,[L_{i},Y_{j}]=-\lambda_i^{j}X_{j}   \\[0.5ex]
        \,[L_{1},L_{2}]=0
    \end{array}$$
    for $i=1,2$ and $j,k=1,\dots,m$,
    }
\end{ex}
\begin{ex}{\rm
    Now consider the metric Lie algebras 
    $\fd(\lambda_1,\lambda_2) := \fd_{\alpha,0}(m,1,2,\lambda)$ for a fixed 
    $\alpha\not=0$ 
    which have
    signature $(2,2+2m+1)$ and were defined in the previous section, too. Here
    $\lambda_1=\lambda(L_1),\ \lambda_2=\lambda(L_2)\in\RR^m$ for a fixed basis $L_1,L_2$ of 
    $\RR^2$ which is choosen in such a way that $\alpha(L_1,L_2)$ is a unit vector. 
    This algebra can be described as follows: $\fd(\lambda_1,\lambda_2)$ 
    for $\lambda_1,\,\lambda_2\in\RR^m$
    is isomorphic to $\RR^{2,2}\oplus\RR^{2m+1}$ as a pseudo-Euclidean vector space, 
    there exists a basis $L_1, L_2, Z_1, Z_2$ of $\RR^{2,2}$ with
    $$\langle L_i,L_j\rangle = \langle Z_i,Z_j\rangle=0,\ \langle Z_i,L_j\rangle
    =\delta_{ij}$$ for $i,j=1,2$ and an orthonormal basis 
    $X_{1},\dots,X_{m},Y_{1},\dots,Y_{m},A_0$ of $\RR^{2m+1}$ such that
    $$
    \begin{array}{l}
        \, Z_{1},Z_{2}\in\fz(\fd(\lambda_1,\lambda_2))  \\[0.5ex]
	\,[A_{0},X_{j}]=[A_{0},Y_{j}]=0\\[0.5ex]
	\,[A_{0},L_{1}]=Z_{2},\quad [A_{0},L_{2}]=-Z_{1}\\[0.5ex]
        \,[X_{j},Y_{k}]=\delta_{jk}(\lambda_1^{j}Z_{1}+
	\lambda_2^{j}Z_{2}),
	\quad [X_{j},X_{k}]=0,\quad [Y_{j},Y_{k}]=0   \\[0.5ex]
        \,[L_{i},X_{j}]=\lambda_i^{j}Y_{j},\quad
        \,[L_{i},Y_{j}]=-\lambda_i^{j}X_{j}   \\[0.5ex]
        \,[L_{1},L_{2}]=A_{0}
    \end{array}$$
    for $i=1,2$ and $j,k=1,\dots,m$      
    }
\end{ex}
\begin{theo}\label{ii}
    Let $(\fg,\ip)$ be an indecomposable metric Lie algebra of
    signature $(2,q)$. If $\fg$ is simple, then $\fg$ is isomorphic
    to $\fsl(2,\RR)$ and $\ip$ is a multiple of the Killing form.
    If $\fg$ is not simple, then the centre $\fz(\fg)$ of $\fg$ is
    one- or two-dimensional and we are in one of the following cases.
    \begin{enumerate}
  \item If $\dim \fz(\fg)=1$, then $q$ is even and $(\fg,\ip)$ is isomorphic to
    $\fd_{L_{2}}(\RR^{1,1},\RR)$ if $q=2$, or to exactly one of the spaces
    $\fd_{L_{2,\lambda}}(\RR^{1,q-1},\RR)$ if $q\ge4$, 
    where
    $$L_2=\left(\begin{array}{cc} 0&1\\1&0 \end{array}\right)\ \mbox{ and }\ 
    L_{2,\lambda}=\left(\begin{array}{ccc} L_2&0&0\\0&0 & -\diag(\lambda)  \\
    0& \diag(\lambda) & 0 \end{array}\right),$$
    $\lambda=(\lambda^1,\lambda^2,\dots,\lambda^r)$ with 
    $0<\lambda^1\le\dots\le\lambda^r$, $r=q/2-1$.
  \item If $\dim \fz(\fg)=2$ and $\dim\fg$ is even, then $q =2m+2$ with $m\ge3$ 
    and $(\fg,\ip)$
    is isomorphic to $\osc(\lambda_1,\lambda_2)$ for some
    $\lambda_1,\lambda_2 \in \RR^{m}$. There is no index $j\in\{1,\dots,m\}$ such
    that $\lambda_1^{j}=\lambda_2^{j}=0$ and
    the set $\{(\lambda^1_i,\lambda_i^2)\mid i=1,\dots,m\}$ is not contained in the 
    union of two 1-dimensional subspaces. 
    Two such Lie algebras $\osc(\lambda_1,\lambda_2)$ and
    $\osc(\bar\lambda_1,\bar\lambda_2)$ are isomorphic if and
    only if
    $$\Span\{\lambda_1,\lambda_2\}=\Span\{\bar\lambda_1,
    \bar\lambda_2\} \quad \mod  \frak
    S_{m}\ltimes(\ZZ_{2})^{m}.$$
  \item If $\dim \fz(\fg)=2$ and $\dim\fg$ is odd, then $q=2m+3$ for $m\ge0$ 
    and $(\fg,\ip)$
    is isomorphic to $\fd(\lambda_1,\lambda_2)$ for some
    $\lambda_1,\lambda_2 \in \RR^{m}$. There is no index 
    $j\in\{1,\dots,m\}$ such
    that $\lambda_1^{j}=\lambda_2^{j}=0$.
    Two such Lie algebras $\fd(\lambda_1,\lambda_2)$ and
    $\fd(\bar\lambda^1,\bar\lambda^2)$ are isomorphic if and
    only if
    $$(\Span\{\lambda_1,\lambda_2\}, \pm\lambda_1\wedge\lambda_2) 
    = (\Span\{\bar\lambda_1, \bar\lambda_2\},\pm\bar \lambda_1\wedge
    \bar\lambda_2) \quad \mod  \frak
    S_{m}\times(\ZZ_{2})^{m}.$$
    \end{enumerate}
    
    Conversely, all metric Lie algebras $\fd_{L_{2}}(\RR^{1,1},\RR)$,   
    $\fd_{L_{2,\lambda}}(\RR^{1,q-1},\RR)$,
    $\osc(\lambda_1,\lambda_2)$, $\fd(\lambda_1,\lambda_2)$ are indecomposable 
    and of signature $(2,q)$ if $\lambda$, $\lambda_1$, $\lambda_2$ satisfy the
    conditions described above. 
\end{theo}
\proof
If $(\fg,\ip)$ is simple, then it is isomorphic to $\fsl(2,\RR)$, because there is no other simple 
Lie algebra whose Killing form has signature $(2,q)$. If $(\fg,\ip)$ is not simple, then by 
Proposition \ref{i2} there exist an abelian Lie
algebra $\fl$,  $\dim\fl=\dim\fz(\fg)\in\{1,2\}$, a semi-Euclidean vector
    space $(\fa, \ip_{\fa})$, an orthogonal representation $\rho$ of
    $\fl$ on $\fa$, and a cocycle $\alpha \in Z^{2}(\fl,\fa)$ satisfying {\rm 
    (\ref{EK})},
    such that $\fd_{\alpha,0}(\fa,\fl,\rho)$ is regular and $(\fg,\ip)$ is 
    isomorphic to
    $\fd_{\alpha,0}(\fa,\fl,\rho)$.

    Let us first consider the case $\dim \fz(\fg)=1$. Then $\fl=\RR$, $(\fa, 
    \ip_{\fa})=\RR^{1,q-1}$ and $\alpha=0$. Hence, $(\fg,\ip)=\fd_A(\RR^{1,q-1},\RR)$ as  
    in Example \ref{ExdA}. Since $\fd_A(\RR^{1,q-1},\RR)$ is regular we have $\ker A=0$ 
    by Lemma \ref{gaga}. By Corollary \ref{ae} regular metric Lie algebras
    $\fd_{A_1}(\RR^{1,q-1},\RR)$ and $\fd_{A_2}(\RR^{1,q-1},\RR)$ are 
    isomorphic if and only if
    \begin{equation}\label{fast} 
    \exists\mu\in\RR^*
    \exists U\in O(1,q-1) : U\,A_1\,U^{-1}=\mu A_2\,.
    \end{equation}
    It is well-known that each equivalence class of bijective maps in $\so(1,q-1)$  
    with respect to relation (\ref{fast}) is represented by exactly one of the maps $L_2$     (if 
    $q=2$) or  
    $L_{2,\lambda}$ with $0<\lambda^1\le\dots\le\lambda^r$, $r=q/2-1$.
    
    If $\dim \fz(\fg)=2$, then the centre is maximal isotropic and the remaining 
    assertions follow 
    directly from Theorem 
    \ref{iso23} . 
\qed

\vspace{0.8cm}
{\footnotesize  
Ines Kath, Martin Olbrich\\
Mathematisches Institut
der Georg-August-Universit\"at G\"ottingen\\
Bunsenstra{\ss}e 3-5, D-37073 G\"ottingen, Germany\\
email: kath@uni-math.gwdg.de, olbrich@uni-math.gwdg.de\\}

\begin{thebibliography}{MMMM}

\bibitem[BK 02]{BK02} H.\,Baum, I.\,Kath, {\sl Doubly Extended Lie Groups -- Curvature, Holonomy 
and Parallel Spinors.} SFB 288 Preprint No. 543, 2002, see also Preprint arXiv:math.DG/0203189.

\bibitem[CP 70]{CP70} M.~Cahen, M.~Parker,  {\sl Sur des classes d'espaces pseudo-riemanniens 
sym\'etriques.} Bull. Soc. Math. Belg. {\bf 22} (1970), 339--354.

\bibitem[CP 80]{CP80} M.~Cahen, M.~Parker, {\sl Pseudo-Riemannian symmetric spaces.} Mem. Amer. 
Math. Soc. {\bf 24} (1980), no. 229.

\bibitem[CW 70]{CW70} M.~Cahen, N.~Wallach, 
{\sl Lorentzian symmetric spaces.} 
Bull. Amer. Math. Soc. {\bf 76} (1970), 585--591.

\bibitem[FS 96]{OF} J.~M.~Figueroa-O'Farrill, S.~Stanciu, {\sl On the structure of symmetric 
self-dual Lie algebras.} J. Math. Phys. {\bf 37} (1996), no. 8, 4121--4134.

\bibitem[Gr 98]{grishkov}
A.~N.~Grishkov,
{\sl Orthogonal modules and nonlinear cohomologies.}
Algebra and Logic {\bf 37} (1998), 294--306.

\bibitem[Gu 64]{gurevich}
G.~B.~Gurevich,
{\sl Foundations of the Theory of Algebraic Invariants.}
P. Noordhoff, Groningen, 1964.

\bibitem[H 00]{hitchin00}
N.~Hitchin,
{\sl The geometry of three-forms in six and seven dimensions.}
Preprint arXiv:math.DG/0010054, 2000.

\bibitem[Kt 92]{katanova}
A.~A.~Katanova,
{\sl Explicit form of certain multivector invariants},
in {\sl Lie Groups, Their Discrete Subgroups and Invariant Theory,}
pages 87--93. Advances in Soviet Mathematices, Vol. 8, AMS Providence, 1992.

\bibitem[M 85]{M85} A.~Medina, {\sl Groupes de Lie munis de m\'etriques bi-invariantes.} Tohoku 
Math. J. (2) {\bf 37} (1985), no. 4, 405--421. 

\bibitem[MR 85]{MR85} A.~Medina, Ph.~Revoy, {\sl Alg\`ebres de Lie et produit scalaire invariant.} 
Ann. Sci. \`Ecole Norm. Sup. (4) {\bf 18} (1985), no. 3, 553--561.

\bibitem[N 02]{N02} Th.~Neukirchner, {\sl Pseudo-Riemannian Symmetric Spaces.} Diplomarbeit, 
Humboldt-Universit\"at zu Berlin, 2002.
\end{thebibliography}
\end{document}